\documentclass[11pt,leqno]{amsart}
\usepackage{amsthm,amsfonts,amssymb,amsmath,oldgerm,mathrsfs,mathabx}
\usepackage[scr=boondoxo]{mathalfa}
\numberwithin{equation}{section}
\usepackage{fullpage,setspace,fancyhdr,bm}
\usepackage{graphicx,psfrag,listings} 
\usepackage{color}

\usepackage[top=30mm,bottom=30mm,left=25mm,right=25mm,a4paper]{geometry}

\usepackage{hyperref}

\newtheorem{definition}{Definition}
\newtheorem{lemma}{Lemma}
\newtheorem{corollary}{Corollary}
\newtheorem{proposition}{Proposition}
\newtheorem{theorem}{Theorem}
\newtheorem{remark}{Remark}

\renewcommand\d{\partial}
\DeclareMathOperator{\beD}{\mathbf{e}}

\DeclareMathOperator{\iD}{i}
\DeclareMathOperator{\eD}{e}
\DeclareMathOperator{\dD}{d}

\def\eps{\varepsilon }

\newcommand{\Id}{\text{Id}}
\newcommand{\I}{\text{I}}
\newcommand{\Zero}{\text{0}}
\newcommand{\One}{\text{\bf 1}}
\DeclareMathOperator{\Tr}{Tr}

\DeclareMathOperator{\sign}{sgn}

\newcommand{\Grho}{G^{\eps,\textrm{ess}}}
\newcommand{\Gtau}{G^{\eps,\textrm{pt}}}
\newcommand{\tGtau}{\tG^{\eps,\textrm{pt}}}
\newcommand{\Gpsi}{G^{\eps,\textrm{p}}}
\newcommand{\Spsi}{s_\eps^{\textrm{p}}}
\newcommand{\Stau}{\tS_\eps}

\newcommand\br{\begin{remark}}
\newcommand\er{\end{remark}}
\newcommand\bp{\begin{pmatrix}}
\newcommand\ep{\end{pmatrix}}
\newcommand{\be}{\begin{equation}}
\newcommand{\ee}{\end{equation}}
\newcommand\ds{\displaystyle}

\newcommand{\bpr}{\begin{proposition}}
\newcommand{\epr}{\end{proposition}}
\newcommand{\bt}{\begin{theorem}}
\newcommand{\et}{\end{theorem}}
\newcommand{\bc}{\begin{corollary}}
\newcommand{\ec}{\end{corollary}}
\newcommand{\bl}{\begin{lemma}}
\newcommand{\el}{\end{lemma}}
\newcommand{\bd}{\begin{definition}}
\newcommand{\ed}{\end{definition}}

\newcommand\R{\mathbf R}
\newcommand\C{\mathbf C}
\newcommand{\N}{\mathbf N}

\newcommand{\RR}{{\mathbb R}}

\newcommand\bA{{\mathbf A}}

\newcommand\bK{{\mathbf K}}
\newcommand\bL{{\mathbf L}}

\newcommand\bR{{\mathbf R}}

\newcommand\bT{{\mathbf T}}

\newcommand\bV{{\mathbf V}}
\newcommand\bW{{\mathbf W}}

\newcommand\bPhi{{\mathbf \Phi}}

\newcommand\bfalpha{{\boldsymbol \alpha}}

\newcommand\uU{{\underline U}}

\newcommand\uu{{\underline u}}

\newcommand\cA{{\mathcal A}}
\newcommand\cB{{\mathcal B}}
\newcommand\cC{{\mathcal C}}
\newcommand\cD{{\mathcal D}}

\newcommand\cL{{\mathcal L}}
\newcommand\cM{{\mathcal M}}
\newcommand\cN{{\mathcal N}}
\newcommand\cO{{\mathcal O}}

\newcommand\ttt{{\widetilde t}}
\newcommand\tu{{\widetilde u}}

\newcommand\tx{{\widetilde x}}

\newcommand\tmu{{\widetilde \mu}}

\newcommand\trho{{\widetilde \rho}}

\newcommand\tsigma{{\widetilde \sigma}}

\newcommand\ttau{{\widetilde \tau}}

\newcommand\tom{{\widetilde \omega}}

\newcommand\tcB{\widetilde{\cB}}
\newcommand\tcN{\widetilde{\cN}}
\newcommand\tuU{\widetilde{\uU}}

\newcommand\tG{\widetilde{G}}

\newcommand\tP{\widetilde{P}}

\newcommand\tS{\widetilde{S}}

\newcommand\tSigma{\widetilde{\Sigma}}

\newcommand\tbA{\widetilde{\bA}}

\newcommand\tbL{\widetilde{\bL}}

\newcommand\tbR{\widetilde{\bR}}

\newcommand\tbV{\widetilde{\bV}}

\newcommand\tbPhi{\widetilde{\bPhi}}

\title{Uniform asymptotic stability\\for convection-reaction-diffusion equations\\in the inviscid limit towards Riemann shocks}

\author{Paul Blochas}
\address{Univ Rennes, CNRS, IRMAR - UMR 6625, F-35000 Rennes, France}
\email{{\tt paul.blochas@univ-rennes1.fr}}
\thanks{Research of P.B. was partially supported by the French region of Brittany.}

\author{L.~Miguel Rodrigues}
\address{
Univ Rennes \& IUF, CNRS, IRMAR - UMR 6625, F-35000 Rennes, France}
\email{{\tt luis-miguel.rodrigues@univ-rennes1.fr}}
\thanks{}

\begin{document}

\maketitle

\begin{abstract}
The present contribution proves the asymptotic orbital stability of viscous regularizations of stable Riemann shocks of scalar balance laws, uniformly with respect to the viscosity/diffusion parameter $\eps$. The uniformity is understood in the sense that all constants involved in the stability statements are uniform and that the corresponding multiscale $\eps$-dependent topology reduces to the classical $W^{1,\infty}$-topology when restricted to functions supported away from the shock location. Main difficulties include that uniformity precludes any use of parabolic regularization to close regularity estimates, that the global-in-time analysis is also spatially multiscale due to the coexistence of nontrivial slow parts with fast shock-layer parts, that the limiting smooth spectral problem (in fast variables) has no spectral gap and that uniformity requires a very precise and unusual design of the phase shift encoding \emph{orbital} stability. In particular, our analysis builds a phase that somehow interpolates between the hyperbolic shock location prescribed by the Rankine-Hugoniot conditions and the non-uniform shift arising merely from phasing out the non-decaying $0$-mode, as in the classical stability analysis for fronts of reaction-diffusion equations.

\vspace{0.5em}


\vspace{0.5em}

\end{abstract}

\tableofcontents

\section{Introduction}

In the present contribution, we prove for the very first time an asymptotic stability result, uniform with respect to the viscosity parameter, for a viscous regularization of a discontinuous traveling-wave of an hyperbolic equation. 

\subsection{The original hyperbolic result}

The purely inviscid result \cite{DR1}, that we extend to the slightly viscous regimes, is itself quite recent. More generally, despite the fact that hyperbolic models are largely used for practical purposes and that for such models singularities such as shocks and characteristic points are ubiquitous, the analysis of nonlinear asymptotic stability of singular traveling-waves of hyperbolic systems is still in its infancy. The state of the art is essentially reduced to a full classification of waves of scalar equations in any dimension \cite{DR1,DR2} (obtained using some significant insights about characteristic points from \cite{JNRYZ}) and the case study of a discontinuous wave without characteristic point for a system of two equations in dimension $1$ \cite{SYZ,YangZumbrun20}. 

Let us stress that, in the foregoing, stability is understood in the sense of Lyapunov, that is, globally in time, and for a topology encoding piecewise smoothness. This is consistent with the fact that concerning stability in the sense of Hadamard, that is, short-time well-posedness, for piecewise-smooth topologies, a quite comprehensive (but not complete) theory is already available even for multidimensional systems; see \cite{Majda1,Majda2,Metivier_cours-chocs,Benzoni-Serre_book}. At this level of regularity, being a weak solution is characterized by a free-interface initial boundary value problem, composed of equations taken in the classical sense in zones of smoothness, and the Rankine-Hugoniot transmission conditions along the free interfaces of discontinuity. 

As is well-known, for hyperbolic equations, weak solutions are not unique and one needs to make an extra choice. The one we are interested in is the most classical one when the extra condition is to be obtained as a vanishing viscosity limit. For scalar equations, in any dimension, since the pioneering work of Kru{\v z}kov \cite{Kruzhkov} (see also \cite[Chapters~4 and~6]{Bressan}), this is known to be sufficient to ensure uniqueness and to be characterized by the so-called entropy conditions, which at our level of smoothness are reduced to inequalities at the free interfaces of discontinuity. For systems, even in dimension~$1$, despite decisive breakthroughs achieved in \cite{Bianchini-Bressan}, such questions are still the object of intensive research; see for instance \cite{Kang-Vasseur}. The present contribution lies at the crossroad of these questions related to the basic definitions of the notion of solution for hyperbolic equations and the ongoing development of a robust general theory for the stability of traveling waves, for which we refer the reader to \cite{Sattinger-book,Henry-geometric,Zumbrun,Sandstede,KapitulaPromislow-stability,JNRZ-conservation}. From the former point of view, the present contribution may be thought as a global-in-time scalar version of \cite{Goodman-Xin,Grenier-Rousset,Rousset}. From the latter point of view, though of a very different technical nature, by many respects, it shares similar goals with other vanishing viscosity stability programs --- see for instance \cite{Bedrossian-Germain-Masmoudi,Herda-Rodrigues} --- and the present contribution is thought as being to \cite{DR1} what \cite{Bedrossian-Masmoudi-Vicol} is to \cite{Bedrossian-Masmoudi}.

We focus on the most basic shock stability result of \cite{DR1}. Consider a scalar balance law in dimension~$1$,
\be\label{eq:hyp}
\d_t u+\d_x(f(u))=g(u)
\ee
with traveling wave solutions $\R\times\R\to\R$, $(t,x)\mapsto \uu(x-(\psi_0+\sigma_0 t))$ with initial shock position $\psi_0\in\R$, speed $\sigma_0\in\R$ and wave profile $\uu$ of Riemann shock type, that is,
\[ \uu(x)=\begin{cases}
\uu_{-\infty} &\text{ if }  x<0\\
\uu_{+\infty} &\text{ if }  x>0
\end{cases}\]
where $(\uu_{-\infty},\uu_{+\infty})\in\R^2$, $\uu_{+\infty}\neq\uu_{-\infty}$. The fact that this does define a weak solution is equivalent to
\begin{align}\label{hyp-weak}
g(\uu_{+\infty})&=0\,,& g(\uu_{-\infty})&=0\,,&
f(\uu_{+\infty})-f(\uu_{-\infty})&=\sigma_0 (\uu_{+\infty}-\uu_{-\infty})\,,
\end{align}
whereas a strict version of entropy conditions may be enforced in Oleinik's form
\begin{equation}\label{hyp-Lax}
\begin{cases}
\qquad\qquad\qquad\sigma_0\,>\,f'(\uu_{+\infty})\,,&\\[0.5em]
\frac{f(\tau\,\uu_{-\infty}+(1-\tau)\,\uu_{+\infty})-f(\uu_{-\infty})}{\tau\,\uu_{-\infty}+(1-\tau)\,\uu_{+\infty}-\uu_{-\infty}}>
\frac{f(\tau\,\uu_{-\infty}+(1-\tau)\,\uu_{+\infty})-f(\uu_{+\infty})}{\tau\,\uu_{-\infty}+(1-\tau)\,\uu_{+\infty}-\uu_{+\infty}}&\qquad\textrm{for any }\ \tau\in(0,1)\,,\\[0.5em]
\qquad\qquad\qquad f'(\uu_{-\infty})\,>\,\sigma_0\,.
\end{cases}
\end{equation}
Requiring a strict version of entropy conditions ensures that they still hold for nearby functions and in particular they disappear at the linearized level. In the foregoing, and throughout the text, for the sake of simplicity, we assume that $f,g\in \cC^\infty(\R)$ though each result only requires a small amount of regularity.

The following statement is one of the alternative versions of \cite[Theorem~2.2]{DR1} described in \cite[Remark~2.3]{DR1}. 

\bt[\cite{DR1}]\label{th:hyp}
Let $(\sigma_0,\uu_{-\infty},\uu_{+\infty})\in\R^3$ define a strictly-entropic Riemann shock of \eqref{eq:hyp} in the above sense. Assume that it is spectrally stable in the sense that
\[
g'(\uu_{+\infty})<0 \quad \text{ and } \quad g'(\uu_{-\infty})<0\,.
\]
There exist $\delta>0$ and $C>0$ such that for any $\psi_0\in\R$ and $v_0\in BUC^1(\R^*)$ satisfying
\[
\|v_0\|_{W^{1,\infty}(\R^*)}\leq\delta\,,
\]
there exists $\psi\in\cC^2(\R^+)$ with initial data $\psi(0)=\psi_0$ such that the entropy solution to~\eqref{eq:hyp}, $u$, generated by the initial data $u(0,\cdot)=(\uu+v_0)(\cdot+\psi_0)$, belongs to $BUC^1(\R_+\times\R\setminus\{\,(t,\psi(t))\,;\,t\geq0\,\})$ and satisfies for any $t\geq 0$
\begin{align*}
\|u(t,\cdot-\psi(t))-\uu\|_{W^{1,\infty}(\R^*)}
+|\psi'(t)-\sigma_0|
&\leq \|v_0\|_{W^{1,\infty}(\R^*)}\,C\, \eD^{\max(\{g'(\uu_{+\infty}),g'(\uu_{-\infty})\})\, t}\,,
\end{align*}
and moreover there exists $\psi_\infty$ such that 
\[
|\psi_\infty-\psi_0|\,\leq \|v_0\|_{L^{\infty}(\R^*)} C\,,
\]
and for any $t\geq 0$
\[
|\psi(t)-\psi_\infty-t\,\sigma_0|\,\leq \|v_0\|_{L^{\infty}(\R^*)} C\, \eD^{\max(\{g'(\uu_{+\infty}),g'(\uu_{-\infty})\})\, t}\,.
\]
\et

In the foregoing, we have used notation $BUC^k(\Omega)$ to denote the set of $\cC^k$ functions over $\Omega$ whose derivatives up to order $k$ are bounded, and uniformly continuous on every connected component of $\Omega$. In other words, $BUC^k(\Omega)$ is the closure of $W^{\infty,\infty}(\Omega)$ for the $W^{k,\infty}(\Omega)$ topology. Working with $BUC^k$ instead of $W^{k,\infty}$ allows to use approximation by smooth functions, an argument ubiquitous in local well-posedness theories, without imposing vanishing at $\infty$.

Note that expressed in classical stability terminology the previous theorem provides asymptotic orbital stability with asymptotic phase. We stress however that the role of phase shifts is here deeper than in the classical stability analysis of smooth waves since it is not only required to provide decay of suitable norms in large-time but also to ensure that these norms are finite locally in time. In particular here there is no freedom, even in finite time, in the definition of phase shifts that need to synchronize discontinuities to allow for comparisons in piecewise smooth topologies.

It is also instructive to consider the corresponding spectral problem. In a moving frame, linearizing from $u(t,x)=\uu(x-(\psi_0+\sigma_0 t)-\psi(t))+v(t,x-(\psi_0+\sigma_0 t)-\psi(t))$ gives a linear IBVP in $(v,\psi)$ 
\begin{align*}
(\d_t&+(f'(\uu_{+\infty})-\sigma_0)\d_x-g'(\uu_{+\infty}))\,v(t,\cdot)\,=\,0\qquad\textrm{on }\RR_+^*\,,\\
(\d_t&+(f'(\uu_{-\infty})-\sigma_0)\d_x-g'(\uu_{-\infty}))\,v(t,\cdot)\,=\,0\qquad\textrm{on }\RR_-^*\,,\\
\psi'(t)
&-\,\left(\frac{f'(\uu_{+\infty})-\sigma_0}{\uu_{+\infty}-\uu_{-\infty}}v(t,0^+)
-\frac{f'(\uu_{-\infty})-\sigma_0}{\uu_{+\infty}-\uu_{-\infty}}v(t,0^-)\right)\,=\,0\,.
\end{align*}
The corresponding spectrum on $BUC^1(\R^*)\times\R$ is 
\[
\left\{\ \lambda\ ;\ \Re(\lambda)\leq \max(\{g'(\uu_{-\infty});g'(\uu_{+\infty})\})\ \right\}\cup\{0\}
\]
and when $\max(\{g'(\uu_{-\infty});g'(\uu_{+\infty})\})<0$, $0$ has multiplicity $1$ (in the sense provided by resolvent singularities) with eigenvector $(0,1)$. This shows that Theorem~\ref{th:hyp} sharply reproduces linear behavior.

\subsection{The vanishing viscosity problem}

Since even the local-in-time notion of solution involves vanishing viscosity approximations, it is natural to wonder whether Theorem~\ref{th:hyp} may have a small-viscosity extension or whether the local-in-time vanishing viscosity limits may be globalized in time about the stable Riemann shocks of Theorem~\ref{th:hyp}. We answer such a question for the following parabolic approximation
\be\label{eq:scalar}
\d_t u+\d_x(f(u))=\eps\,\d_x^2 u+g(u)\,.
\ee 
Note that solutions to \eqref{eq:scalar} are smooth (not uniformly in $\eps$) so that techniques based on free-interfaces IBVP formulations for \eqref{eq:hyp} cannot easily be extended to the study of \eqref{eq:scalar}. In the reverse direction, to gain a better control on smoothness of solutions to \eqref{eq:scalar}, it is expedient to introduce fast variables 
\[
u(t,x)=\tu\underbrace{\left(\frac{t}{\eps},\frac{x}{\eps}\right)}_{(\ttt,\tx)}
\]
that turn \eqref{eq:scalar} into
\be\label{eq:scalar-fast}
\d_\ttt \tu+\d_\tx(f(\tu))=\d_\tx^2\tu+\eps\,g(\tu)\,.
\ee
We stress however that this is indeed in original variables $(t,x)$ that we aim at proving a uniform result. In particular, a large part of the analysis is focused on distinctions between norms that get large and norms that get small when going from slow to fast variables. For a closely related discussion we refer the reader to \cite{Kang-Vasseur_bis,Kang-Vasseur}.

In order to carry out the extension, the first step is to elucidate the existence of traveling waves to \eqref{eq:scalar} near $\uu$. A preliminary observation in this direction is that the formal $\eps\to0$ limit of \eqref{eq:scalar-fast} does possess a smooth traveling-wave solution $(\ttt,\tx)\mapsto \uU_0(\tx-\sigma_0\,\ttt)$ of speed $\sigma_0$ and profile $\uU_0$ such that
\begin{align*}
\lim_{-\infty}\uU_0&=\uu_{-\infty}\,,&
\lim_{+\infty}\uU_0&=\uu_{+\infty}\,,
\end{align*}
simply obtained by solving
\begin{align*}
\uU_0(0)&=\frac{\uu_{-\infty}+\uu_{+\infty}}{2}\,,&
\uU_0'&\,=\,f(\uU_0)-f(\uu_{+\infty})
-\sigma_0\,(\uU_0-\uu_{+\infty})\,.
\end{align*}
We recall that $\sigma_0$ is tuned to ensure $f(\uu_{-\infty})-f(\uu_{+\infty})
-\sigma_0\,(\uu_{-\infty}-\uu_{+\infty})=0$ and observe that the Oleinik's entropy conditions imply that $\uU_0$ is strictly monotonous. This $\eps=0$ viscous profile is often called viscous shock layer and plays the role of a short-time free-interface boundary layer. This simple limiting fast profile may be perturbed to yield profiles for \eqref{eq:scalar-fast} hence for \eqref{eq:scalar}. To state such a perturbation result with optimal spatial decay rates, we introduce, for $\eps\geq0$,
\begin{align*}
\theta_\eps^r
&:=\frac{1}{2}|f'(\uu_{+\infty})-\sigma_\eps|
+\frac{1}{2}\sqrt{(f'(\uu_{+\infty})-\sigma_\eps)^2+4\,\eps\,|g'(\uu_{+\infty})|}\,,\\
\theta_\eps^\ell
&:=\frac{1}{2}\,(f'(\uu_{-\infty})-\sigma_\eps)
+\frac{1}{2}\sqrt{(f'(\uu_{-\infty})-\sigma_\eps)^2+4\,\eps\,|g'(\uu_{-\infty})|}\,.
\end{align*}

\bpr\label{p:profile}
Under the assumptions of Theorem~\ref{th:hyp}, for any $0<\alpha^\ell<\theta_0^\ell$, $0<\alpha^r<\theta_0^r$ and $k_0\in\N^*$, there exist $\eps_0>0$ and $C_0>0$ such that there exist a unique $(\sigma_\eps,\uU_\eps)$, with $\uU_\eps\in\cC^2(\R)$,
\begin{align*}
\uU_\eps(0)&=\frac{\uu_{-\infty}+\uu_{+\infty}}{2}\,,&
(f(\uU_\eps)-\sigma_\eps\,\uU_\eps)'&=\uU_\eps''+\eps\,g(\uU_\eps)\,,
\end{align*}
and
\begin{align*}
|\sigma_\eps-\sigma_0|
+\|\eD^{\alpha^\ell|\,\cdot\,|}(\uU_\eps-\uU_0)\|_{W^{1,\infty}(\R_-)}
+\|\eD^{\alpha^r\,\cdot\,}(\uU_\eps-\uU_0)\|_{W^{1,\infty}(\R_+)}
&\leq C_0\,\eps\,,
\end{align*}
and, moreover, there also holds
\begin{align*}
\|\eD^{\theta_\eps^\ell|\,\cdot\,|}(\uU_\eps-\uu_{-\infty})
-\eD^{\theta_0^\ell|\,\cdot\,|}(\uU_0-\uu_{-\infty})\|_{L^{\infty}(\R_-)}
&\leq\, C_0\,\eps\,,\\
\|\eD^{\theta_\eps^r\,\cdot\,}(\uU_\eps-\uu_{+\infty})
-\eD^{\theta_0^r\,\cdot\,}(\uU_0-\uu_{+\infty})\|_{L^{\infty}(\R_+)}
&\leq\, C_0\,\eps\,,\\
\|\eD^{\theta_\eps^\ell|\,\cdot\,|}\uU_\eps^{(k)}
-\eD^{\theta_0^\ell|\,\cdot\,|}\uU_0^{(k)}\|_{L^{\infty}(\R_-)}
&\leq\, C_0\,\eps\,,&1\leq k\leq k_0\,,\\
\|\eD^{\theta_\eps^r\,\cdot\,}\uU_\eps^{(k)}
-\eD^{\theta_0^r\,\cdot\,}\uU_0^{(k)}\|_{L^{\infty}(\R_+)}
&\leq\, C_0\,\eps\,,&1\leq k\leq k_0\,.
\end{align*}
\epr

Note that a traveling-wave $(t,x)\mapsto \uu_\eps(x-(\psi_0+\sigma_\eps t))$ with $\psi_0\in\R$ arbitrary, is obtained from $\uU_\eps$ through
\[
\uu_\eps(x)\,:=\,\uU_\eps\left(\frac{x}{\eps}\right)
\]
and that, uniformly in $\eps$,
\begin{align*}
|\uu_\eps(x)-\uu(x)|&\lesssim \eD^{-\theta_\eps^\ell\,\frac{|x|}{\eps}}\,,& x<0\,,\\
|\uu_\eps(x)-\uu(x)|&\lesssim \eD^{-\theta_\eps^r\,\frac{x}{\eps}}\,,& x>0\,,\\
|\uu^{(k)}_\eps(x)|&\lesssim \frac{1}{\eps^k}\eD^{-\theta_\eps^\ell\,\frac{|x|}{\eps}}\,,& x<0\,,&\ k\geq 1\,,\\
|\uu^{(k)}_\eps(x)|&\lesssim \frac{1}{\eps^k}\eD^{-\theta_\eps^r\,\frac{x}{\eps}}\,,& x>0\,,&\ k\geq 1\,.\\
\end{align*}

We prove Proposition~\ref{p:profile} in Appendix~\ref{s:profiles}. The existence and uniqueness part with suboptimal spatial rates follows from a rather standard Lyapunov-Schmidt argument. We stress however that it is crucial for our linear and nonlinear stability analyses to gain control on $\uU'_\eps$ with sharp spatial decay rates. We obtain the claimed upgrade from suboptimal to optimal rates essentially as a corollary to the refined spectral analysis needed to carry out the nonlinear study. We point out that, despite the fact that the literature on the subject is quite extensive --- see for instance \cite{Harterich_hyperbolic,Harterich_canard,Crooks-Mascia,Crooks,Gilding} and references therein ---, we have not found there an existence result with the level of generality needed here, that is, including non-convex fluxes and yielding optimal spatial decay rates. 

With the existence of $\eps$-versions of traveling waves in hands, the next natural question is whether these are spectrally stable. It is settled by standard arguments, as expounded in \cite{KapitulaPromislow-stability}, combining direct computations of the essential spectrum with Sturm-Liouville theory. The latter uses crucially that $\uU'_\eps$ is monotonous, a consequence of the Oleinik's entropy conditions. The upshot is that, in slow original variables, the spectrum of the linearization about $\uu_\eps$ in a co-moving frame, acting on $BUC^1(\R)$, is stable and exhibits a spectral gap between the simple eigenvalue $0$ and the rest of the spectrum of size $\min(\{|g'(\uu_{+\infty})|,|g'(\uu_{-\infty})|\})+\cO(\eps)$. Note that in fast variables the spectral gap is of size $\eps\times\min(\{|g'(\uu_{+\infty})|,|g'(\uu_{-\infty})|\})+\cO(\eps^2)$ Details of the latter are given in Section~\ref{s:spectral}.

The real challenge is \emph{uniform} nonlinear asymptotic stability. Indeed, if one removes the uniformity requirement, nonlinear stability follows from spectral stability by now well-known classical arguments as expounded in \cite{Sattinger-book,Henry-geometric,Sandstede,KapitulaPromislow-stability}, and initially developed in, among others, \cite{Sattinger_one,Sattinger_two,Henry-geometric,Kapitula,Wu-Xing,Xing}. Since the limit is singular, it is worth spelling out what we mean by uniform stability. There are two closely related parts in the requirement. Explicitly, on initial data, 
\begin{enumerate}
\item the most obvious one is that the restriction on the sizes of allowed initial perturbations (encoded by the smallness of $\delta$ in Theorem~\ref{th:hyp}) should be uniform with respect to $\eps$, so that the lower bound on the size of the basin of attraction provided by the analysis is nontrivial in the limit $\eps\to0$;
\item the second one is more intricate\footnote{But our result satisfies a much simpler and stronger version of the requirement.}, it states that the $\eps$-dependent norms, say $\|\cdot\|_{(\eps)}$, used to measure this smallness (in slow original variables) should be controlled by an $\eps$-independent norm for functions supported away from the shock, so that in particular for any $v\in \cC^\infty_c(\R)$ supported in $\R^*$, $\limsup_{\eps\to0}\|v\|_{(\eps)}<+\infty$.
\end{enumerate}
On the control of solutions arising from perturbations, we impose similar constraints but with upper bounds replacing lower bounds in the requirements. Constraints on the control of solutions ensure that the bounds provide a nontrivial control whereas constraints on the control of initial data ensure that nontrivial perturbations are allowed. 

It may be intuitive that the stringer the norm is the larger the size of the basin of attraction is since a qualitatively better control is offered by the topology. In the present case, the discussion is on the amount of localization encoded by the norm since, though this is somewhat hidden, time decay is controlled by initial spatial localization (as opposed to cases where regularity drives decay as for instance in \cite{Bedrossian-Masmoudi,Bedrossian-Masmoudi-Vicol,Bedrossian-Germain-Masmoudi}). To offer a quantitative insight, let us use as in \cite{Herda-Rodrigues} a simple ODE as a toy model to predict the size constraints. Consider the stability of $y\equiv 0$ for $y'=-\tau\,y+\rho\,y^2$ where $\tau>0$ encodes the size of the spectral gap and $\rho>0$ measures the size of nonlinear forcing. For such an equation, a ball of radius $r_0$ and center $0$ is uniformly attracted to $0$ provided that $r_0<\tau/\rho$. Now, if one considers \eqref{eq:scalar} directly in $BUC^1(\R)$ (or any reasonable unweighted topology) and forgets about issues related to phase definitions and possible regularity losses, the spectral gap offered by a linearization about $\uu_\eps$ is of order $1$ whereas the forcing by nonlinear terms is of order $\eps^{-1}$ (since this is the size of $\uu_\eps'$) 
hence the rough prediction of a basin of size $\cO(\eps)$. Yet, working with weights such as $\eD^{-\theta\,\eps^{-\alpha}\,|x|}$, for some sufficiently small $\theta>0$ and some $0<\alpha\leq 1$, moves the spectrum to increase the size of the gap to the order $\eps^{-\alpha}$ yielding the expectation of an $\cO(\eps^{1-\alpha})$ basin. Note that the choice $\alpha=1$ would provide a uniform size and is consistent with the size of viscous shock layers but it would force initial perturbations to be located in an $\cO(\eps)$ spatial neighborhood of the shock location. 

The foregoing simple discussion predicts quite accurately\footnote{Actually it is even a bit optimistic for the unweighted and $\alpha<1$ cases.} what could be obtained by applying the most classical parabolic strategy to the problem at hand. The failure of the classical strategy may also be read on the deeply related, but not equivalent, fact that it uses the phase only to pull out the contribution of nonlinear terms through the spectrally non-decaying $0$-mode. This is inconsistent with the stronger role of the phase for the hyperbolic problem, all the spectrum contributing to the phase in the latter case. A completely different approach is needed.

Additional strong signs of the very challenging nature of the uniform stability problem may also be gathered from the examination of the viscous layer stability problem, that is, the stability of $(\ttt,\tx)\mapsto\uU_0(\tx-\sigma_0\,\ttt)$ as a solution to \eqref{eq:scalar-fast} with $\eps=0$. The problem has been extensively studied, see for instance \cite{Liu,Goodman,Goodman_bis,Jones-Gardner-Kapitula,Kreiss-Kreiss,Howard_lin,Howard_nonlin} for a few key contributions and \cite{Zumbrun} for a thorough account. The spectrum of the linearization includes essential spectrum touching the imaginary axis at $0$, which is still an eigenvalue, so that the decay is not exponential but algebraic and requires a trade-off, localization against decay, as for the heat equation. This is a consequence of a conservative nature of the equation, but the conservative structure may also be used to tame some of the apparent difficulties. To give one concrete example: one may remove the embedded eigenvalue $0$ from the essential spectrum by using the classical antiderivative trick, dating back at least to \cite{Matsumura-Nishihara}, either directly at the nonlinear level under the restriction of zero-mean perturbations (as in \cite{Goodman_ter} or in \cite{Matsumura-Nishihara} for a system case) or only to facilitate the linear analysis as in \cite{Howard_lin,Howard_nonlin}. In fast variables, turning on $\eps>0$ moves the essential to the left, creating an $\cO(\eps)$ spectral gap but breaks the conservative structure thus rendering almost impossible, and at least quite inconvenient, the use of classical conservative tools. Our stability analysis requires a description as detailed as the one of \cite{Howard_lin,Howard_nonlin} and, without the antiderivative trick at hand, this involves the full machinery of \cite{Zumbrun-Howard,Zumbrun-Howard_erratum}. Roughly speaking, one of the main outcomes of our detailed spectral analysis, expressed in fast variables, is that the $\eps$-proximity of essential spectrum and $0$-eigenvalue induces that the essential spectrum has an impact of size $1/\eps$ on the linear time-evolution, but that at leading-order the algebraic structure of the essential-spectrum contribution is such that it may be absorbed in a suitably designed phase modulation. Note that this is consistent with the fact that, in fast variables, variations in shock positions are expected to be of size $1/\eps$ and with the fact that, in slow variables, the phase is involved in the resolution of all the hyperbolic spectral problems, not only the $0$-mode.

To summarize and extend the discussion so far, we may hope 
\begin{enumerate}
\item to overcome the discrepancy between the Rankine-Hugoniot prescription of the phase and the pure $0$-mode modulation, and to phase out the hidden singularity caused by the proximity of essential spectrum and $0$ eigenvalue, by carefully identifying the most singular contribution of the essential spectrum as phase variations and including this in a carefully designed phase; 
\item to guarantee uniform nonlinear decay estimates provided that we can ensure that, in sow variables, nonlinear terms of size $1/\eps$ also come with a spectral-gap enhancing factor $\eD^{-\theta\,|x|/\eps}$ (for some $\theta>0$).
\end{enumerate}
The latter expectation is motivated by the fact that it is indeed the case for terms forced by $\uu'_\eps$ 
but we need to prove that it is so also for stiff terms caused by the derivatives of the perturbation itself. Concerning the latter, we stress that even if one starts with a very gentle perturbation supported away from the shock the nonlinear coupling instantaneously creates stiff parts of shock-layer type in the perturbation thus making it effectively multi-scale. 

There remains a somewhat hidden issue, that we have not discussed so far. Along the foregoing discussion, we have done as if we could use Duhamel principle based on a straight-forward linearization, as in classical semilinear parabolic problems. Yet, here, closing nonlinear estimates in regularity by using parabolic regularization either explicitly through gains of derivatives or indirectly through $L^q\to L^p$, $q<p$, mapping properties, effectively induces losses in power of $\eps$ in an already $\eps$-critical problem thus is completely forbidden. Instead, we estimate
\begin{itemize}
\item the variation in shock position $\psi$, the shape variation $v$ and the restriction of its derivative $\d_xv$ to an $\cO(\eps)$ neighborhood of the shock location through Duhamel formula and linear decay estimates;
\item the remaining part of $\d_xv$ by a suitably modified Goodman-type hyperbolic energy estimate.
\end{itemize}
The latter energy estimate is similar in spirit to those in \cite{Goodman_ter,Rodrigues-Zumbrun,YangZumbrun20} but the hard part of its design is precisely in going from a classical hyperbolic estimate that would work in the complement of an $\cO(1)$ neighborhood of the shock location to a finely tuned estimate covering the complement of an $\cO(\eps)$ neighborhood, since this is required for the combination with a lossless parabolic regularization argument. Moreover, there are two more twists in the argument: on one hand we need the estimate to include weights encoding the multi-scale nature of $\d_xv$ ; on the other hand, for the sake of sharpness, to remain at the $\cC^1$ level of regularity, we actually apply the energy estimates on a suitable nonlinear version of $\d_xv$ so that they yield $L^\infty$ bounds for $\d_xv$. 

The arguments sketched above, appropriately worked out, provide the main result of the present paper. To state such results, we introduce multi-scale weights and corresponding norms: for $k\in\N$, $\eps>0$, and $\theta\geq0$, 
\begin{align*}
\omega_{k,\eps,\theta}(x)
&:=
\frac{1}{1+\frac{1}{\eps^k}\,\eD^{-\theta\,\frac{|x|}{\eps}}}\,,&
\|v\|_{W_{\eps,\theta}^{k,\infty}(\R)}
&=\sum_{j=0}^k\,\|\omega_{j,\eps,\theta}\,\d_x^jv\|_{L^{\infty}(\R)}
\end{align*}
Note that
\begin{enumerate}
\item Each norm $\|\cdot\|_{W_{\eps,\theta}^{k,\infty}(\R)}$ is equivalent to any standard norm on $\|\cdot\|_{W^{k,\infty}(\R)}$ but non uniformly in $\eps$ and that the uniformity is restored if one restricts it to functions supported in the complement of a fixed neighborhood of the origin.
\item The norm $\|\cdot\|_{W_{\eps,\theta}^{0,\infty}(\R)}$ is uniformly equivalent to $\|\cdot\|_{L^{\infty}(\R)}$.
\item If $\theta<\min(\{\theta_0^\ell,\theta_0^r\})$ then $\|\uu_\eps-\uu\|_{W_{\eps,\theta}^{k,\infty}(\R)}$ is bounded uniformly with respect to $\eps$.
\end{enumerate}

\bt\label{th:main-sl}
Enforce the assumptions and notation of Theorem~\ref{th:hyp} and Proposition~\ref{p:profile}.\\
There exists $\theta_0>0$ such that for any $0<\theta\leq\theta_0$, there exist $\eps_0>0$, $\delta>0$ and $C>0$ such that for any $0<\eps\leq\eps_0$, any $\psi_0\in\R$ and any $v_0\in BUC^1(\R)$ satisfying
\[
\|v_0\|_{W_{\eps,\theta}^{1,\infty}(\R)}
\leq\delta\,,
\]
there exists $\psi\in\cC^1(\R^+)$ with initial data $\psi(0)=\psi_0$ such that the strong\footnote{We ensure $u\in BUC^{0}(\R_+;BUC^{1}(\R))\cap \cC^\infty(\R_+^*;BUC^{\infty}(\R))$.} solution to~\eqref{eq:scalar}, $u$, generated by the initial data $u(0,\cdot)=(\uu_\eps+v_0)(\cdot+\psi_0)$, is global in time and satisfies for any $t\geq 0$
\begin{align*}
\|u(t,\cdot-\psi(t))-\uu_\eps\|_{W_{\eps,\theta}^{1,\infty}(\R)}
+|\psi'(t)-\sigma_\eps|
&\leq \|v_0\|_{W_{\eps,\theta}^{1,\infty}(\R)}\,C\, \eD^{\max(\{g'(\uu_{+\infty}),g'(\uu_{-\infty})\})\, t}\,,
\end{align*}
and moreover there exists $\psi_\infty$ such that 
\[
|\psi_\infty-\psi_0|\,\leq  \|v_0\|_{W_{\eps,\theta}^{1,\infty}(\R)} C\,,
\]
and for any $t\geq 0$
\[
|\psi(t)-\psi_\infty-t\,\sigma_\eps|\,\leq  \|v_0\|_{W_{\eps,\theta}^{1,\infty}(\R)} C\, \eD^{\max(\{g'(\uu_{+\infty}),g'(\uu_{-\infty})\})\, t}\,.
\]
\et

Among the many variations and extensions of Theorem~\ref{th:hyp} provided in \cite{DR1}, the simplest one to extend to a uniform small viscosity result is \cite[Proposition~2.5]{DR1} that proves that the exponential time decay also holds for higher order derivatives without further restriction on sizes of perturbations. It does not require any new insight besides the ones used to prove Theorem~\ref{th:main-sl} and we leave it aside only to cut unnecessary technicalities.

Likewise, one may obtain in an even more direct way, that is, up to immaterial changes, exponential damping of norms encoding further slow spatial localization. To give an  explicit example, let us extend notation $W_{\eps,\theta}^{k,\infty}$, $L^\infty$, into $W_{\eps,\theta,\theta'}^{k,\infty}$, $L^\infty_{\theta'}$, accordingly to weights 
\begin{align*}
\omega_{k,\eps,\theta,\theta'}(x)
&:=
\frac{1}{\eD^{-\theta'\,|x|}+\frac{1}{\eps^k}\,\eD^{-\theta\,\frac{|x|}{\eps}}}\,,&
\omega_{\theta'}(x)
&:=\eD^{\theta'\,|x|}\,,&
\end{align*}
with $\theta'\geq0$ arbitrary. One may prove for instance that for any $\theta'\geq 0$ there exist $C_{\theta'}$ and $\eps_{\theta'}>0$ such that, under the sole further restrictions $0<\eps\leq \eps_{\theta'}$ and $\eD^{\theta'\,|\,\cdot\,|}\,v_0\in L^\infty(\R)$, there holds
\begin{align*}
\|u(t,\cdot-\psi(t))-\uu_\eps\|_{L_{\theta'}^{\infty}(\R)}
&\leq \|v_0\|_{L_{\theta'}^{\infty}(\R)}\,C_\theta\, \eD^{\max(\{g'(\uu_{+\infty}),g'(\uu_{-\infty})\})\, t}\,.
\end{align*}
One point in considering these weighted topologies is that, when $\theta'>0$, $L^\infty_{\theta'}$ is continuously embedded in $L^1\cap L^\infty$, so that an estimate on $\|u(t,\cdot-\psi(t))-\uu\|_{L^p}$ is provided by the combination of the foregoing bound with the already known bound
\[
\|\uu_\eps-\uu\|_{L^p(\R)}\lesssim \eps^{\frac1p}\,.
\]

\subsection{Outline and perspectives}

The most natural nontrivial extensions of Theorems~\ref{th:hyp}/\ref{th:main-sl} that we have chosen to leave for future work concern on one hand the parabolic regularization by quasilinear terms and on the other hand planar Riemann shocks in higher spatial dimensions (see \cite[Theorem~3.4]{DR1} for the hyperbolic case). We expect many parts of the present analysis to be directly relevant in quasilinear or multiD cases but we also believe that their treatments would also require sufficiently many new arguments to deserve a separate treatment.

In the multidimensional case, even the outcome is expected to be significantly different. In this direction, let us point out that the hyperbolic spectral problem is critical in the stronger sense that the spectrum includes the whole imaginary axis, instead of having an intersection with the imaginary axis reduced to $\{0\}$. This may be tracked back to the fact that the linearized Rankine-Hugoniot equation takes the form of a transport equation in transverse variables for the phase. Consistently, as proved in \cite[Theorem~3.4]{DR1}, for the hyperbolic problem, perturbing a planar shock may lead asymptotically in large time to another non-planar Riemann shock sharing the same constant-states. This may still be interpreted as a space-modulated asymptotic stability result, in the sense coined in \cite{JNRZ-conservation} and thoroughly discussed in \cite{R,R_Roscoff,R_linKdV,DR2}. A similar phenomenon is analyzed for scalar \emph{conservation} laws in \cite{Serre_scalar}.

Concerning the quasilinear case, the main new difficulty is expected to arise from the fact that, to close the argument, one needs to prove that the $L^\infty$ decay of $\eps\,\d_x^2v$, where $v$ still denotes the shape variation, is at least as good as the one of $\d_xv$. A priori, outside the shock layer this leaves the freedom to pick some initial typical size $\eps^{-\eta_0}$, $\eta_0\in [0,1]$, for $\d_x^2v$ and to try to propagate it. Indeed, roughly speaking, in the complement of an $\cO(\eps)$ neighborhood of the shock location, this $L^\infty$ propagation stems from arguments similar to the ones sketched above for $\d_xv$. The key difference is that now one cannot complete it with a bound obtained through Duhamel formula since this would involve an $L^\infty$ bound on $\d_x^3v$. Thus the quasilinear study seems to require to be able to close an estimate for $\d_x^2v$ entirely with energy-type arguments, a highly non-trivial task. 


In another direction, we expect that the study of waves with characteristic points, as arising in the full classification obtained in \cite{DR2} for scalar balance laws, should not only involve some new patches here and there but follow very different routes and thus will require significantly new insights even at a general abstract level. As a strong token of this expectation, we point out that regularity is expected to play a paramount role there since, at the hyperbolic level, the regularity class chosen deeply modifies the spectrum when a characteristic point is present in the wave profile ; see \cite{JNRYZ,DR2}.

\smallskip

The rest of the paper is organized as follows. We have decided to shift the derivation of wave profile asymptotics, proving Proposition~\ref{p:profile}, to Appendix~\ref{s:profiles}, because we believe that the backbone of the paper is stability and provide it mostly for completeness' sake. Next section contains a detailed examination of the required spectral preliminaries. The following one explains how to use these to obtain a practical representation of the linearized time-evolution. Though we mostly follow there the arguments in \cite{Zumbrun-Howard}, with some twists here and there, we provide a detailed exposition for two distinct reasons. The first one is that we need to track in constructions which parts are $\eps$-uniform and which parts are not, a crucial point in our analysis. The second one is that most of the papers of the field requiring a detailed analysis, as we do, are either extremely long \cite{Zumbrun-Howard} or cut in a few long pieces \cite{Mascia-Zumbrun_hyp-par,Mascia-Zumbrun_hyp-par_bis} and we want to save the reader from back-and-forth consultations of the literature. This makes our analysis essentially self-contained (up to basic knowledge of spectral analysis) and we believe that it could serve as a gentle introduction to the latter massive literature. Note however, that, to keep the paper within a reasonable size, we only expound the bare minimum required by our analysis. After these two preliminary sections, we enter into the technical core of the paper, with first a section devoted to detailed linear estimates, including the identification of most-singular parts of the time-evolution as phase variations, and then a section devoted to nonlinear analysis, including adapted nonlinear maximum principles proved through energy estimates and the proof of Theorem~\ref{th:main-sl}. 

\section{Spectral analysis}\label{s:spectral}




We investigate stability for traveling waves introduced in Proposition~\ref{p:profile}. We have chosen to carry out all our proofs within co-moving fast variables. Explicitly, we introduce new unknowns and variables through\footnote{Note the slight co-moving inconsistency with the introduction.}  
\[
u(t,x)=\tu\underbrace{\left(\frac{t}{\eps},\frac{x-\sigma_\eps\,t}{\eps}\right)}_{(\ttt,\tx)}.
\]
However, since we never go back to the original slow variables, we drop tildes on fast quantities from now on. One reason to opt for the fast variables is that it provides a simpler reading of size dependencies on $\eps$.

Therefore our starting point is
\be\label{eq:scalar-sl}
\d_t u+\d_x(f(u)-\sigma_\eps\,u)=\d_x^2 u+\eps\,g(u)\,,
\ee
about the stationary solution $\uU_\eps$. Accordingly we consider the operator
\be\label{def:operator}
\cL_\eps:=-\d_x((f'(\uU_\eps)-\sigma_\eps)\,\cdot\,)+\d_x^2+\eps g'(\uU_\eps)
\ee
on $BUC^0(\R)$ with domain $BUC^2(\R)$.

Though the elements we provide are sufficient to reconstruct the classical theory, the reader may benefit from consulting \cite{KapitulaPromislow-stability} for background on spectral analysis specialized to nonlinear wave stability. In particular, we shall make extensive implicit use of the characterizations of essential spectrum in terms of endstates of wave profiles and of the spectrum at the right-hand side\footnote{We picture the complex plane with the real axis pointing to the right and the imaginary axis pointing to the top.} of the essential spectrum\footnote{There are (at least) two reasonable definitions of essential spectrum, either through failure of satisfying Fredholm property or through failure of satisfying Fredholm property with zero index. In the context of semigroup generators both definitions provide the same right-hand boundary thus the conventional choice is immaterial to stability issues.} in terms of zeroes of Evans' functions. The reader is referred to \cite{Kato,Davies} for less specialized, basic background on spectral theory.

The backbone of the theory is the interpretation of spectral properties of one-dimensional differential operators in terms of spatial dynamics and a key-part of the corresponding studies is the investigation of exponential dichotomies. It starts with the identification between the eigenvalue equation 
\[
(\lambda-\cL_\eps)\,v\,=\,0
\]
and the system of ODEs
\[
\frac{\dD}{\dD x}\bV(x)\,=\,\bA_\eps(\lambda,x)\,\bV(x)
\]
for the vector\footnote{The use of flux variables is not necessary but it simplifies a few computations here and there.} $\bV=(v,\d_xv-(f'(\uU_\eps)-\sigma_\eps)\,v)$ where
\be\label{def:Aeps}
\bA_\eps(\lambda,x)
\,:=\,\bp
f'(\uU_\eps)-\sigma_\eps&1\\
\lambda-\eps\,g'(\uU_\eps)&0
\ep\,.
\ee
For later use, we shall denote $\bPhi_\eps^\lambda(x,y)$ the corresponding solution operators, mapping datum at point $y$ to value at point $x$.

The essential spectrum is characterized in terms of matrices $\bA_\eps^r(\lambda):=\bA^\eps(\lambda;\uu_{+\infty})$ and $\bA_\eps^\ell(\lambda):=\bA^\eps(\lambda;\uu_{-\infty})$ with
\begin{align}\label{def:Aeps_rl}
\bA^\eps(\lambda;u)
&\,:=\,\bp
f'(u)-\sigma_\eps&1\\
\lambda-\eps\,g'(u)&0
\ep\,.
\end{align}
Eigenvalues of $\bA^\eps(\lambda;u)$ are given by
\be\label{def:mupm}
\mu_{\pm}^\eps(\lambda;u):=\frac{f'(u)-\sigma_\eps}{2}
\pm\sqrt{\frac{(f'(u)-\sigma_\eps)^2}{4}+\lambda-\eps\,g'(u)}
\ee
and are distinct when $\lambda\neq\eps\,g'(u)-\tfrac14\,(f'(u)-\sigma_\eps)^2$. In this case, the matrix may be diagonalized as 
\[
\bA^\eps(\lambda;u)\,=\,\bp \bR_+^\eps(\lambda;u)&\bR_-^\eps(\lambda;u)\ep\,
\bp \mu_+^\eps(\lambda;u)&0\\0&\mu_-^\eps(\lambda;u) \ep\,
\bp \bL_+^\eps(\lambda;u)\\\bL_-^\eps(\lambda;u)\ep
\]
with
\begin{align}\label{def:Vpm}
\bR_{\pm}^\eps(\lambda;u)&:=\bp 1\\-\mu_{\mp}^\eps(\lambda;u)\ep\,,&
\bL_{\pm}^\eps(\lambda;u)&:=\frac{\bp \pm\mu_{\pm}^\eps(\lambda;u)&\pm 1\ep}{
\mu_+^\eps(\lambda;u)-\mu_-^\eps(\lambda;u)}\,.
\end{align}
The eigenvalues $\mu_{\pm}^\eps(\lambda;u)$ have distinct real parts when $\lambda$ does not belong to
\[
\cD_\eps(u):=\eps\,g'(u)-\tfrac14\,(f'(u)-\sigma_\eps)^2+\R^-\,.
\]
All our spectral studies will take place far from the half-lines $\cD_\eps(\uu_{+\infty})\cup\cD_\eps(\uu_{-\infty})$, that correspond to the set termed \emph{absolute spectrum} in \cite{KapitulaPromislow-stability}.

From now on, throughout the text, we shall use $\sqrt{\cdot}$ to denote the determination of the square root on $\C\setminus\R^-$ with positive real part.

\subsection{Conjugation to constant coefficients}

Our starting point is a conjugation of spectral problems to a piecewise constant coefficient spectral problem. This is mostly relevant in compact zones of the spectral plane and in the literature by Kevin Zumbrun and his collaborators this is known as a \emph{gap lemma} --- since a gap or in other words an exponential dichotomy is the key assumption ---; see for instance \cite[Lemma~2.6]{MZ_AMS} for a version relevant for the present analysis. Since we need to ensure uniformity in $\eps$ for the case at hand we provide both a statement and a proof.

\bpr\label{p:gap-lemma}
Let $K$ be a compact subset of $\C\setminus\cD_0(\uu_{+\infty})$. There exist positive constants $(\eps_0,C,\theta)$ such that there exists\footnote{As follows from the proof, $P_\eps^r(\lambda,\cdot)$ is defined as soon as $\lambda\notin\cD_\eps(\uu_{+\infty})\cup\cD_\eps(\uu_{-\infty})$.} a smooth map
\[
P^r\,:\,[0,\eps_0]\times K\times \R\mapsto GL_2(\C)\,,\qquad
(\eps,\lambda,x)\mapsto P_\eps^r(\lambda,x)
\] 
locally uniformly analytic in $\lambda$ on a neighborhood of $K$ and such that, for any $(\eps,\lambda,x)\in [0,\eps_0]\times K\times [0,+\infty)$,
\begin{align*}
\|P_\eps^r(\lambda,x)-\I_2\|&\leq C\,\eD^{-\theta\,|x|}\,,&
\|(P_\eps^r(\lambda,x))^{-1}-\I_2\|&\leq C\,\eD^{-\theta\,|x|}\,,
\end{align*}
and, for any $(\eps,\lambda,x,y)\in [0,\eps_0]\times K\times (\R_+)^2$,
\[
\bPhi_\eps^\lambda(x,y)\,=\,
P_\eps^r(\lambda,x)\,\eD^{(x-y)\,\bA_\eps^r(\lambda)}\,
(P_\eps^r(\lambda,y))^{-1}\,.
\]
\epr

The same argument applies to the conjugation on $(-\infty,0]$ with the flow of $\bA_\eps^\ell(\lambda)$ and defines a conjugation map denoted $P^\ell$ from now on.

\begin{proof}
The proof is essentially a quantitative "cheap" gap lemma --- conjugating only one trajectory instead of solution operators --- but applied in $\cM_2(\C)$ instead of $\C^2$.

Let us first observe that it is sufficient to define $P^r$ on $[0,\eps_0]\times K\times [x_0,+\infty)$ for some suitably large $x_0$. Indeed then one may extend $P^r$ by
\[
P_\eps^r(\lambda,x)\,:=\,
\bPhi_\eps^\lambda(x,x_0)\,P_\eps^r(\lambda,x_0)\,\eD^{(x_0-x)\,\bA_\eps^r(\lambda)}\,
\]
and bounds are extended by a continuity-compactness argument. Likewise the uniformity in $\eps$ is simply derived from a continuity-compactness argument since the construction below is continuous at the limit $\eps=0$. Note moreover that in the large-$x$ regime the bound on $(P_\eps^r(\lambda,x))^{-1}$ may be derived from the bound on $P_\eps^r(\lambda,x)$ by using properties of the inverse map. 

The requirements on $P_\eps^r(\lambda,\cdot)$ are equivalent to the fact that it converges exponentially fast to $\I_2$ at $+\infty$ (uniformly in $\eps$) and that it satisfies for any $x$
\[
\frac{\dD}{\dD x}P_\eps^r(\lambda,x)
\,=\,\cA_\eps^r(\lambda)(P_\eps^r(\lambda,x))
+(\bA_\eps(\lambda,x)-\bA_\eps^r(\lambda))\,P_\eps^r(\lambda,x)
\]  
where $\cA_\eps^r(\lambda):=\cA^\eps(\lambda;\uu_{+\infty})$ is a linear operator on $\cM_2(\C)$ defined through
\[
\cA^\eps(\lambda;u)(P):=[\bA^\eps(\lambda;u),P]
=\bA^\eps(\lambda;u)P-P\bA^\eps(\lambda;u)\,.
\]
When $\lambda\notin\cD_\eps(u)$, $\cA^\eps(\lambda;u)$ admits 
\[
(\bR_+^\eps(\lambda;u)\bL_-^\eps(\lambda;u),\quad
\bR_-^\eps(\lambda;u)\bL_+^\eps(\lambda;u),\quad
\bR_+^\eps(\lambda;u)\bL_+^\eps(\lambda;u),\quad
\bR_-^\eps(\lambda;u)\bL_-^\eps(\lambda;u))
\]
as a basis of eigenvectors corresponding to eigenvalues
\[
(\mu_+^\eps(\lambda;u)-\mu_-^\eps(\lambda;u),\quad
-(\mu_+^\eps(\lambda;u)-\mu_-^\eps(\lambda;u)),\quad
0,\quad
0)\,.
\]
Note that $\I_2$ always lies in the kernel of $\cA^\eps(\lambda;u)$. We denote by $\Pi_u^\eps(\lambda;u)$, $\Pi_s^\eps(\lambda;u)$, $\Pi_0^\eps(\lambda;u)$ the corresponding spectral projections respectively on the unstable space, the stable space and the kernel of $\cA^\eps(\lambda;u)$ and for further later study we point out that they are given as
\begin{align*}
\Pi_u^\eps(\lambda;u)(P)&=\bL_+^\eps(\lambda;u)\,P\,\bR_-^\eps(\lambda;u)\quad\bR_+^\eps(\lambda;u)\bL_-^\eps(\lambda;u)\,,\\
\Pi_s^\eps(\lambda;u)(P)&=\bL_-^\eps(\lambda;u)\,P\,\bR_+^\eps(\lambda;u)\quad\bR_-^\eps(\lambda;u)\bL_+^\eps(\lambda;u)\,,\\
\Pi_0^\eps(\lambda;u)(P)&=\bL_+^\eps(\lambda;u)\,P\,\bR_+^\eps(\lambda;u)\quad\bR_+^\eps(\lambda;u)\bL_+^\eps(\lambda;u)\,\\
&\ +\bL_-^\eps(\lambda;u)\,P\,\bR_-^\eps(\lambda;u)\quad\bR_-^\eps(\lambda;u)\bL_-^\eps(\lambda;u)\,.
\end{align*}

Then the result follows when $\theta$ is sufficiently small and $x_0$ is sufficiently large from a use of the implicit function theorem on
\begin{align*}
P_\eps^r(\lambda,x)
&\,=\,\I_2\,-\int_x^{+\infty}\Pi_0^\eps(\lambda;\uu_{+\infty})\left(
(\bA_\eps(\lambda,y)-\bA_\eps^r(\lambda))\,P_\eps^r(\lambda,y)\right)\,\dD y\\
&-\int_x^{+\infty}\eD^{-(y-x)\,(\mu_+^\eps(\lambda;\uu_{+\infty})-\mu_-^\eps(\lambda;\uu_{+\infty}))}\,\Pi_u^\eps(\lambda;\uu_{+\infty})\left(
(\bA_\eps(\lambda,y)-\bA_\eps^r(\lambda))\,P_\eps^r(\lambda,y)\right)\,\dD y\\
&+\int_{x_0}^x\eD^{-(x-y)\,(\mu_+^\eps(\lambda;\uu_{+\infty})-\mu_-^\eps(\lambda;\uu_{+\infty}))}\,\Pi_s^\eps(\lambda;\uu_{+\infty})\left(
(\bA_\eps(\lambda,y)-\bA_\eps^r(\lambda))\,P_\eps^r(\lambda,y)\right)\,\dD y
\end{align*}
with norm control on matrix-valued maps through 
\[
\sup_{x\geq x_0}\eD^{\theta\,|x|}\|P(x)-\I_2\|\,.
\] 
\end{proof}

\br
The properties of the foregoing proposition do not determine $P^r$ uniquely. The normalizing choice made in the proof is $\Pi_s^\eps(\lambda;\uu_{+\infty})\left(P_\eps^r(\lambda,x_0)\right)=\Zero_2$ but we could have replaced $\Zero_2$ with any analytic choice of an element of the stable space of $\cA^\eps(\lambda;\uu_{+\infty})$.
\er

\br
The proposition is sufficient to prove classical results about the determination of the essential spectrum from endstates spectra.
\er

We now investigate possible failure of uniformity in the regime of large spectral parameters. In the literature by Kevin Zumbrun and his collaborators, similar purposes are achieved through comparison of unstable manifolds with their frozen-coefficients approximations by a type of lemma termed there \emph{tracking lemma}; see for instance \cite{Humpherys-Lyng-Zumbrun,BJRZ}. The  rationale is that to large-frequencies smooth coefficients seem almost constant and thus may be treated in some adiabatic way, a fact ubiquitous in high-frequency/semiclassical analysis. 

We follow here a different path and rather effectively build a conjugation as in the foregoing gap lemma. The first step is a suitable scaling to ensure some form of uniformity in the large-$x$ contraction argument of the proof of Proposition~\ref{p:gap-lemma}. The second-step is a high-frequency approximate diagonalization combined with an explicit solving of the leading-order part of the system ensuring that in the large-frequency regime the latter construction could actually be carried out with $x_0=0$. 
\bpr\label{p:tracking-lemma}
There exist positive constants $(\eps_0,C,\theta,\delta)$ such that setting 
\[
\Omega_\delta:=\left\{\,\lambda\,;\ 
\Re\left(\sqrt{\lambda}\right)
\,\geq \frac{1}{\delta}
\right\}
\]
there exists a smooth map
\[
P^{r,HF}\,:\,[0,\eps_0]\times \Omega_\delta\times \R\mapsto GL_2(\C)\,,\qquad
(\eps,\lambda,x)\mapsto P_\eps^{r,HF}(\lambda,x)
\] 
locally uniformly analytic in $\lambda$ on a neighborhood of $\Omega_\delta$ and such that, for any $(\eps,\lambda,x)\in [0,\eps_0]\times \Omega_\delta\times [0,+\infty)$,
\begin{align*}
\left\|\bp 1&0\\0&\frac{1}{\sqrt{\lambda}}\ep\,P_\eps^{r,HF}(\lambda,x)\,\bp 1&0\\0&\sqrt{\lambda}\ep\,
-\eD^{-\frac12\int_x^{+\infty}(f'(\uU_\eps(y))-f'(\uu_{+\infty}))\,\dD y}\,\I_2
\right\|
&\leq C\,\frac{\eD^{-\theta\,|x|}}{\Re\left(\sqrt{\lambda}\right)}\,,\\
\left\|\bp 1&0\\0&\frac{1}{\sqrt{\lambda}}\ep
(P_\eps^{r,HF}(\lambda,x))^{-1}
\bp 1&0\\0&\sqrt{\lambda}\ep
-\eD^{\frac12\int_x^{+\infty}(f'(\uU_\eps(y))-f'(\uu_{+\infty}))\,\dD y}\,\I_2\right\|
&\leq C\,\frac{\eD^{-\theta\,|x|}}{\Re\left(\sqrt{\lambda}\right)}\,,
\end{align*}
and, for any $(\eps,\lambda,x,y)\in [0,\eps_0]\times \Omega_\delta\times \R^2$,
\[
\bPhi_\eps^\lambda(x,y)\,=\,
P_\eps^{r,HF}(\lambda,x)\,\eD^{(x-y)\,\bA_\eps^r(\lambda)}\,
(P_\eps^{r,HF}(\lambda,y))^{-1}\,.
\]
\epr

We point out that for our main purposes we do not need to identify explicitly the leading order part of the conjugation. 

As for Proposition~\ref{p:gap-lemma} the same argument applies to the conjugation on $(-\infty,0]$ with the flow of $\bA_\eps^\ell(\lambda)$ and defines a conjugation map denoted $P^{\ell,HF}$ from now on.

\begin{proof}
As a preliminary remark, we observe that the condition $\Re\left(\sqrt{\lambda}\right)\gg1$ and $\Re(\mu_+^\eps(\lambda;\uu_{+\infty})-\mu_-^\eps(\lambda;\uu_{+\infty}))\gg 1$ are equivalent, with uniform control of one by the other and vice versa.

Scaling $P_\eps^r$ to $Q_\eps^r$ defined as 
\[
Q_\eps^r(\lambda,\cdot)
:=\bp 1&0\\0&\mu_+^\eps(\lambda;\uu_{+\infty})-\mu_-^\eps(\lambda;\uu_{+\infty})\ep^{-1}\,
P_\eps^r(\lambda,\cdot)
\,\bp 1&0\\0&\mu_+^\eps(\lambda;\uu_{+\infty})-\mu_-^\eps(\lambda;\uu_{+\infty})\ep
\]
removes high-frequency singularities by replacing $\bR_{\pm}^\eps(\lambda;\uu_{+\infty})$ and $\bL_{\pm}^\eps(\lambda;\uu_{+\infty})$ with
\begin{align*}
\bp 1\\-\frac{\mu_{\mp}^\eps(\lambda;\uu_{+\infty})}{
\mu_+^\eps(\lambda;\uu_{+\infty})-\mu_-^\eps(\lambda;\uu_{+\infty})}\ep\,,&&
\bp \frac{\pm\mu_{\pm}^\eps(\lambda;\uu_{+\infty})}{
\mu_+^\eps(\lambda;\uu_{+\infty})-\mu_-^\eps(\lambda;\uu_{+\infty})}&\pm 1\ep\,,
\end{align*}
whereas the only other effect is the replacement of $\bA_\eps(\lambda,y)-\bA_\eps^r(\lambda)$ with
\[
\bp
f'(\uU_\eps)-f'(\uu_{+\infty})&0\\[0.25em]
-\eps\,\frac{g'(\uU_\eps)-g'(\uu_{+\infty})}{
\mu_+^\eps(\lambda;\uu_{+\infty})-\mu_-^\eps(\lambda;\uu_{+\infty})}&0
\ep\,.
\]

At this stage, let us choose coordinates to identify $\cM_2(\C)$ with $\C^4$ in such a way that $\C^2\times\{0\}^2$, $\{0\}\times\C\times\{0\}$, $\{0\}^3\times\C$, $(1,0,0,0)$ and $(0,1,0,0)$ correspond respectively --- after scaling and choice of coordinates --- to the kernel of $\cA_\eps^r(\lambda)$, its unstable space, its stable space, $\I_2$ and $\dfrac{1}{\mu_+^\eps(\lambda;\uu_{+\infty})-\mu_-^\eps(\lambda;\uu_{+\infty})}A^r_{\eps ,HF}$. Then the problem to be solved takes the form
\begin{align*}
\frac{\dD}{\dD x}\bfalpha(x)&=\cB_0(x)\,\bfalpha(x)
+\cO(\eD^{-\theta\,|x|})\,\bp \beta\\\gamma\ep(x)\\
\frac{\dD}{\dD x}\beta(x)&=
(\mu_+^\eps(\lambda;\uu_{+\infty})-\mu_-^\eps(\lambda;\uu_{+\infty})
+\omega_+(x))\,\beta(x)
+\cO(\eD^{-\theta\,|x|})\,\bp \bfalpha\\\gamma\ep(x)\\
\frac{\dD}{\dD x}\gamma(x)&=
-(\mu_+^\eps(\lambda;\uu_{+\infty})-\mu_-^\eps(\lambda;\uu_{+\infty})
+\omega_-(x))\,\gamma(x)
+\cO(\eD^{-\theta\,|x|})\,\bp \bfalpha\\\beta\ep(x)
\end{align*}
with $(1,0,0,0)$ as limiting value at $+\infty$, for some\footnote{Along the proof we allow ourselves to change the precise value of $\theta$ from line to line.} $\theta>0$, where $\cB_0(x)$, $\omega_+(x)$, $\omega_-(x)$ are also of the form $\cO(\eD^{-\theta\,|x|})$. 

It follows that when $\mu_+^\eps(\lambda;\uu_{+\infty})-\mu_-^\eps(\lambda;\uu_{+\infty})$ is sufficiently large, by a further change of variables differing from $\I_4$ by a block off-diagonal term
\[
\bp \bfalpha_{bis}\\\beta_{bis}\\\gamma_{bis}\ep(x)
\,=\,\left(\I_4+\cO\left(\frac{\eD^{-\theta\,|x|}}{\mu_+^\eps(\lambda;\uu_{+\infty})-\mu_-^\eps(\lambda;\uu_{+\infty})}\right)\right)\bp \bfalpha\\\beta\\\gamma\ep(x)
\]
one may transform the problem to
\begin{align*}
\frac{\dD}{\dD x}\bfalpha_{bis}(x)&=\tcB_0(x)\,\bfalpha_{bis}(x)
+\cO\left(\frac{\eD^{-\theta\,|x|}}{\mu_+^\eps(\lambda;\uu_{+\infty})-\mu_-^\eps(\lambda;\uu_{+\infty})}\right)\,\bp \beta_{bis}\\\gamma_{bis}\ep(x)\\
\frac{\dD}{\dD x}\beta_{bis}(x)&=
(\mu_+^\eps(\lambda;\uu_{+\infty})-\mu_-^\eps(\lambda;\uu_{+\infty})
+\tom_{+}(x))\,\beta_{bis}(x)\\
&\hspace{7em}
+\cO\left(\frac{\eD^{-\theta\,|x|}}{\mu_+^\eps(\lambda;\uu_{+\infty})-\mu_-^\eps(\lambda;\uu_{+\infty})}\right)\,\bp \bfalpha_{bis}\\\gamma_{bis}\ep(x)\\
\frac{\dD}{\dD x}\gamma_{bis}(x)&=
-(\mu_+^\eps(\lambda;\uu_{+\infty})-\mu_-^\eps(\lambda;\uu_{+\infty})
+\tom_{-}(x))\,\gamma_{bis}(x)\\
&\hspace{7em}
+\cO\left(\frac{\eD^{-\theta\,|x|}}{\mu_+^\eps(\lambda;\uu_{+\infty})-\mu_-^\eps(\lambda;\uu_{+\infty})}\right)\,\bp \bfalpha_{bis}\\\beta_{bis}\ep(x)
\end{align*}
with $(1,0,0,0)$ as limiting value at $+\infty$, where $\tcB_0-\cB_0$, $\tom_{-}-\omega_{-}$ and $\tom_{+}-\omega_{+}$ are all of the form
\[
\cO\left(\frac{\eD^{-\theta\,|x|}}{\mu_+^\eps(\lambda;\uu_{+\infty})-\mu_-^\eps(\lambda;\uu_{+\infty})}\right)\,.
\]

Now, we point out that there is a single solution to the leading-order part
\[
\frac{\dD}{\dD x}\bfalpha_{main}(x)=\tcB_0(x)\,\bfalpha_{main}(x)
 \]
 with $(1,0)$ as limiting value at $+\infty$. This follows from a fixed point argument on $x\geq x_0$ for $x_0$ large followed by a continuation argument. We may be even more explicit. Indeed an explicit computation yields 
\[
\tcB_0(x)\,=\,\cB_{main}(x)+\cO\left(\frac{\eD^{-\theta\,|x|}}{\mu_+^\eps(\lambda;\uu_{+\infty})-\mu_-^\eps(\lambda;\uu_{+\infty})}\right)
\]
with 
\[
\cB_{main}(x)\,\bp 1\\0\ep\,=\,\frac12(f'(\uU_\eps(x))-f'(\uu_{+\infty}))\,\bp 1\\0\ep
\]
so that
\[
\bfalpha_{main}(x)\,=\,\eD^{-\frac12\int_x^{+\infty}(f'(\uU_\eps(y))-f'(\uu_{+\infty}))\,\dD y}\,\bp 1\\0\ep
+\cO\left(\frac{\eD^{-\theta\,|x|}}{\mu_+^\eps(\lambda;\uu_{+\infty})-\mu_-^\eps(\lambda;\uu_{+\infty})}\right)\,.
\] 

The proof is thus achieved by a fixed point argument on a problem of type
\begin{align*}
\bfalpha_{bis}(x)
&\,=\,\bfalpha_{main}(x)\,-\int_x^{+\infty}\Phi_0(x,y)\,
\cO\left(\frac{\eD^{-\theta\,|y|}}{\mu_+^\eps(\lambda;\uu_{+\infty})-\mu_-^\eps(\lambda;\uu_{+\infty})}\right)\,\bp \beta_{bis}\\\gamma_{bis}\ep(y)\,\dD y\\
\beta_{bis}(x)&=
-\int_x^{+\infty}\eD^{-(y-x)\,(\mu_+^\eps(\lambda;\uu_{+\infty})-\mu_-^\eps(\lambda;\uu_{+\infty}))+\int_y^x\tom_+}\,
\cO\left(\frac{\eD^{-\theta\,|y|}}{\mu_+^\eps(\lambda;\uu_{+\infty})-\mu_-^\eps(\lambda;\uu_{+\infty})}\right)\,\bp \bfalpha_{bis}\\\gamma_{bis}\ep(y)\,\dD y\\
\gamma_{bis}(x)&
=\int_{0}^x\eD^{-(x-y)\,(\mu_+^\eps(\lambda;\uu_{+\infty})-\mu_-^\eps(\lambda;\uu_{+\infty}))-\int_y^x\tom_-}\,
\cO\left(\frac{\eD^{-\theta\,|y|}}{\mu_+^\eps(\lambda;\uu_{+\infty})-\mu_-^\eps(\lambda;\uu_{+\infty})}\right)\,\bp \bfalpha_{bis}\\\beta_{bis}\ep(y)\,\dD y
\end{align*}
where $\Phi_0$ denotes the solution operator associated with $\tcB_0$.
\end{proof}

In the following we shall complete Proposition~\ref{p:tracking-lemma} that provides $P^{r,HF}$ on $\Omega_\delta$ with an application of Proposition~\ref{p:gap-lemma} on 
\[
K_\delta:=
\left\{\,\lambda\,;\ 
d(\lambda,\cD_0(\uu_{+\infty}))\,\geq\,\delta
\quad\textrm{and}\quad
|\lambda|\leq \frac{1}{\delta}
\right\}\,.
\]
When $\delta$ is sufficiently small $\Omega_\delta$ and $K_\delta$ overlaps. Yet a priori $P^r$ and $P^{r,HF}$ differ from each other even in regions where both exist. Fortunately the implied possible mismatch disappears at the level of Green functions.

\subsection{Evans' function and its asymptotics}

Now wherever it makes sense we set
\begin{align*}
\bV^{r,s}_\eps(\lambda,x)&:=\eD^{x\,\mu_-^\eps(\lambda;\uu_{+\infty})}\,P^r_\eps(\lambda,x)\,\bR_-^\eps(\lambda;\uu_{+\infty})\,,\\
\bV^{r,u}_\eps(\lambda,x)&:=\eD^{x\,\mu_+^\eps(\lambda;\uu_{+\infty})}\,P^r_\eps(\lambda,x)\,\bR_+^\eps(\lambda;\uu_{+\infty})\,,\\
\bV^{\ell,s}_\eps(\lambda,x)&:=\eD^{x\,\mu_-^\eps(\lambda;\uu_{-\infty})}\,P^\ell_\eps(\lambda,x)\,\bR_-^\eps(\lambda;\uu_{-\infty})\,,\\
\bV^{\ell,u}_\eps(\lambda,x)&:=\eD^{x\,\mu_+^\eps(\lambda;\uu_{-\infty})}\,P^\ell_\eps(\lambda,x)\,\bR_+^\eps(\lambda;\uu_{-\infty})\,,
\end{align*}
and similarly for $\bV^{r,s,HF}_\eps$, $\bV^{r,u,HF}_\eps$, $\bV^{\ell,s,HF}_\eps$ and $\bV^{\ell,u,HF}_\eps$. Note that notation is used here to recall stable and unstable spaces, we also use them in areas of the spectral plane where they do not match with stable and unstable spaces. Instead, this fits analytic continuation of generators of stable/unstable spaces.

Correspondingly we define the Evans' function
\be\label{def:Evans}
D_\eps(\lambda):=\det\bp\bV^{r,s}_\eps(\lambda,0)&\bV^{\ell,u}_\eps(\lambda,0)\ep
\ee
and its high-frequency counterpart $D^{HF}_\eps$. Note that we define the Evans function at point $0$ but on one hand, we do not make use of any particular property due to normalization so that the point $0$ could be replaced with any other point and on the other hand relations between Evans functions at different points are simply derived from Liouville's formula for Wronskians. For instance,
\be\label{eq:Wronskian}
\det\bp\bV^{r,s}_\eps(\lambda,x)&\bV^{\ell,u}_\eps(\lambda,x)\ep
\,=\,D_\eps(\lambda)\,\eD^{\int_0^x\Tr(\bA_\eps(\lambda,\cdot))}
\,=\,D_\eps(\lambda)\,\eD^{\int_0^x(f'(\uU_\eps)-\sigma_\eps)}\,.
\ee

A simple corollary to Proposition~\ref{p:tracking-lemma} is

\bc\label{c:Evans-large}
Uniformly in $\eps$ (sufficiently small)
\[
\lim_{\Re(\sqrt{\lambda})\to\infty}\,\frac{D^{HF}_\eps(\lambda)}{\sqrt{\lambda}}\,=\,
2\,\eD^{-\frac12\int_0^{+\infty}(f'(\uU_\eps(y))-f'(\uu_{+\infty}))\,\dD y
+\frac12\int_{-\infty}^0(f'(\uU_\eps(y))-f'(\uu_{-\infty}))\,\dD y}\,.
\]
\ec

To complete Corollary~\ref{c:Evans-large}, we derive information on compacts sets of $\lambda$ in the limit $\eps\to0$.

\bpr\label{p:Sturm-Liouville}
There exists $\eta_0>0$ such that for any $\delta>0$ there exist positive $(\eps_0,c_0)$ such that for any $\eps\in[0,\eps_0]$, $D_\eps(\cdot)$ is well-defined on
\[
K_{\eta_0,\delta}
:=
\left\{\,\lambda\,;\ 
d(\lambda,(-\infty,-\eta_0])\,\geq\,\min\left(\left\{\delta,\frac{\eta_0}{2}\right\}\right)
\quad\textrm{and}\quad|\lambda|\leq \frac{1}{\delta}
\right\}\,
\]
and for any $\lambda\in K_{\eta_0,\delta}$,
\begin{align*}
|D_\eps(\lambda)|\,\geq\,c_0\,\min(\{1,|\lambda|\})\,.
\end{align*}
\epr

\begin{proof}
We derive the result from Sturm-Liouville theory and regularity in $\eps$. To apply Sturm-Liouville theory, we introduce the weight
\[
\omega_\eps(x):=\eD^{\frac12\int_0^x(f'(\uU_\eps(y))-\sigma_\eps)\,\dD y}\,.
\]
We observe that considered as an operator on $L^2(\R)$ with domain $H^2(\R)$ the operator 
\[
L_\eps:=\frac{1}{\omega_\eps}\cL_\eps\left(\omega_\eps\,\cdot\,\right)
\]
is self-adjoint and in the region of interest it possesses no essential spectrum and its eigenvalues agree in location and algebraic multiplicity with the roots of $D_\eps$. As a consequence the zeroes of $D_\eps$ are real and since $\uU'_\eps/\omega_\eps$ is a nowhere-vanishing eigenvector for the eigenvalue $0$, $0$ is a simple root of $D_\eps$ and $D_\eps$ does not vanish on $(0,+\infty)$. From here the corresponding bound is deduced through a continuity-compactness argument in $\eps$.
\end{proof}

\section{Green functions}\label{s:Green}

Now we use the introduced spectral objects to obtain representation formulas for linearized solution operators.

\subsection{Duality}

To begin with, to provide explicit formulas for spectral Green functions to be introduced below, we extend to dual problems the conclusions of Section~\ref{s:spectral}. Note that the duality we are referring to is not related to any particular choice of a specific Banach space but rather distributional/algebraic.

To begin with, we introduce the formal adjoint
\be\label{def:dual-operator}
\cL_\eps^{adj}:=(f'(\uU_\eps)-\sigma_\eps)\d_x+\d_x^2+\eps g'(\uU_\eps)
\ee
and note that for any sufficiently smooth $v$, $w$, and any points $(x_0,x_1)$
\be\label{eq:duality}
\bV\cdot
\bp 0&-1\\1&0\ep
\bp w\\\d_x w\ep\,(x_1)
-\bV\cdot
\bp 0&-1\\1&0\ep
\bp w\\\d_x w\ep\,(x_0)
\,=\,\int_{x_0}^{x_1}(w\,\cL_\eps v-v\,\cL_\eps^{adj}w)\,
\ee
with $\bV=(v,\d_xv-(f'(\uU_\eps)-\sigma_\eps)\,v)$. As a first simple consequence of \eqref{eq:duality} note that if $(\lambda,y,v_y^r,v_y^\ell)$ are such that $(\lambda-\cL_\eps)v_y^r=0$ and $(\lambda-\cL_\eps)v_y^\ell=0$, then the function
\[
\varphi_y\,:\ \R\to\C\,,\qquad
x\mapsto \begin{cases}
v_y^r(x)&\quad \textrm{if }x>y\\
v_y^\ell(x)&\quad \textrm{if }x<y
\end{cases}
\]
solves $(\lambda-\cL_\eps)\varphi_y\,=\,\delta_y$ if and only if
\begin{align*}
v_y^r(y)&\,=\,v_y^\ell(y)\,,\\
\d_xv_y^r(y)-(f'(\uU_\eps)(y)-\sigma_\eps)\,v_y^r(y)\,
&=\d_xv_y^\ell(y)-(f'(\uU_\eps)(y)-\sigma_\eps)\,v_y^\ell(y)+1\,.
\end{align*}
Specializing to the tensorized case where $v_y^r(x)=v^r(x)\,\alpha(y)$, $v_y^\ell(x)=v^\ell(x)\,\beta(y)$, note that the foregoing conditions are equivalent to
\begin{align*}
\bp \bV^r(y)&\bV^\ell(y)\ep\,
\bp
\alpha(y)\\
-\beta(y)
\ep\,=\,\bp 0\\1\ep
\end{align*}
where $\bV^{\sharp}=(v^{\sharp},\d_xv^{\sharp}-(f'(\uU_\eps)-\sigma_\eps)\,v^{\sharp})$, $\sharp\in\{r,\ell\}$.
Hence, we need to find vectors satisfying some orthogonality property to identify the inverse of the matrix : \begin{align*}
\bp \bV^r(y)&\bV^\ell(y)\ep\,
\end{align*}
To go further, we identify 
\[
(\lambda-\cL_\eps^{adj})\,w\,=\,0
\]
and the system of ODEs
\[
\frac{\dD}{\dD x}\bW(x)\,=\,\tbA_\eps(\lambda,x)\,\bW(x)
\]
for the vector $\bW=(w,\d_x w)$ where
\be\label{def:tAeps}
\tbA_\eps(\lambda,x)
\,:=\,\bp
0&1\\
\lambda-\eps\,g'(\uU_\eps)&
-(f'(\uU_\eps)-\sigma_\eps)
\ep\,.
\ee
Note that 
\[
\tbA_\eps(\lambda,x)\,=\,\bA_\eps(\lambda,x)\,-\,(f'(\uU_\eps)-\sigma_\eps)\,\I_2
\]
so that all the proofs of Section~\ref{s:spectral} purely based on limiting-matrices spectral gaps  arguments apply equally well to the corresponding dual problems under the exact same assumptions. Alternatively one may derive results on dual problems by using directly the relation between solution operators
\[
\tbPhi_\eps^\lambda(x,y)\,=\,\bPhi_\eps^\lambda(x,y)\,\eD^{-\int_y^x(f'(\uU_\eps)-\sigma_\eps)}\,.
\]
Here and elsewhere throughout the text from now on we denote with a $\,\widetilde{\,\,}\,$ all quantities arising from dual problems. Let us point out that our choices lead to
\begin{align*}
\tmu_{\pm}^\eps(\lambda;u)&
\,=\,\mu_{\pm}^\eps(\lambda;u)-\,(f'(u)-\sigma_\eps)\,,&
\tbR_{\pm}^\eps(\lambda;u)&\,=\,\bR_{\pm}^\eps(\lambda;u)\,,&
\tbL_{\pm}^\eps(\lambda;u)&\,=\,\bL_{\pm}^\eps(\lambda;u)\,,
\end{align*}
and
\begin{align*}
\tP^r_\eps(\lambda,x)&
\,=\,P^r_\eps(\lambda,x)\,\eD^{\int_x^{+\infty}(f'(\uU_\eps)-f'(\uu_{+\infty}))}\,,&
\tP^\ell_\eps(\lambda,x)&
\,=\,P^\ell_\eps(\lambda,x)\,\eD^{-\int_{-\infty}^x(f'(\uU_\eps)-f'(\uu_{-\infty}))}\,,&
\end{align*}
(and likewise for high-frequency versions).

\bpr\label{p:duality-mean}
Let $K$ be a compact subset of $\C\setminus(\cD_0(\uu_{+\infty})\cup\cD_0(\uu_{+\infty}))$. There exists $\eps_0>0$ such that there exist smooth maps
\begin{align*}
\tau^r\,&:\,[0,\eps_0]\times K\mapsto\C\,,\ 
(\eps,\lambda)\mapsto \tau_\eps^r(\lambda)\,,&
\rho^r\,&:\,[0,\eps_0]\times K\mapsto\C\,,\ 
(\eps,\lambda)\mapsto \rho_\eps^r(\lambda)\,,\\
\tau^\ell\,&:\,[0,\eps_0]\times K\mapsto\C\,,\ 
(\eps,\lambda)\mapsto \tau_\eps^\ell(\lambda)\,,&
\rho^\ell\,&:\,[0,\eps_0]\times K\mapsto\C\,,\ 
(\eps,\lambda)\mapsto \rho_\eps^\ell(\lambda)\,,\\
\ttau^r\,&:\,[0,\eps_0]\times K\mapsto\C\,,\ 
(\eps,\lambda)\mapsto \ttau_\eps^r(\lambda)\,,&
\trho^r\,&:\,[0,\eps_0]\times K\mapsto\C\,,\ 
(\eps,\lambda)\mapsto \trho_\eps^r(\lambda)\,,\\
\ttau^\ell\,&:\,[0,\eps_0]\times K\mapsto\C\,,\ 
(\eps,\lambda)\mapsto \ttau_\eps^\ell(\lambda)\,,&
\trho^\ell\,&:\,[0,\eps_0]\times K\mapsto\C\,,\ 
(\eps,\lambda)\mapsto \trho_\eps^\ell(\lambda)\,,
\end{align*}
locally uniformly analytic in $\lambda$ on a neighborhood of $K$ and such that, for any $(\eps,\lambda)\in [0,\eps_0]\times K$, for any $x\in\R$
\begin{align*}
\bV_\eps^{r,s}(\lambda,x)&=\,\rho^r_\eps(\lambda)\,\bV_\eps^{\ell,s}(\lambda,x)
+\tau^r_\eps(\lambda)\,\bV_\eps^{\ell,u}(\lambda,x)\,,\\
\bV_\eps^{\ell,u}(\lambda,x)&=\,\rho^\ell_\eps(\lambda)\,\bV_\eps^{r,u}(\lambda,x)
+\tau^\ell_\eps(\lambda)\,\bV_\eps^{r,s}(\lambda,x)\,,\\
\tbV_\eps^{r,s}(\lambda,x)&=\,\trho^r_\eps(\lambda)\,\tbV_\eps^{\ell,s}(\lambda,x)
+\ttau^r_\eps(\lambda)\,\tbV_\eps^{\ell,u}(\lambda,x)\,,\\
\tbV_\eps^{\ell,u}(\lambda,x)&=\,\trho^\ell_\eps(\lambda)\,\tbV_\eps^{r,u}(\lambda,x)
+\ttau^\ell_\eps(\lambda)\,\tbV_\eps^{r,s}(\lambda,x)\,,\\
\end{align*}
and
\begin{align*}
\rho^r_\eps(\lambda)&\,=\,
D_\eps(\lambda)\,\frac{\eD^{-\int_{-\infty}^0(f'(\uU_\eps)-f'(\uu_{-\infty}))}}{\mu_+^\eps(\lambda;\uu_{-\infty})-\mu_-^\eps(\lambda;\uu_{-\infty})}\,,\\
\rho^\ell_\eps(\lambda)&\,=\,
D_\eps(\lambda)\,\frac{\eD^{\int_0^{+\infty}(f'(\uU_\eps)-f'(\uu_{+\infty}))}}{\mu_+^\eps(\lambda;\uu_{+\infty})-\mu_-^\eps(\lambda;\uu_{+\infty})}\,,\\
\trho^r_\eps(\lambda)&\,=\,
D_\eps(\lambda)\,\frac{\eD^{\int_0^{+\infty}(f'(\uU_\eps)-f'(\uu_{+\infty}))}}{\mu_+^\eps(\lambda;\uu_{-\infty})-\mu_-^\eps(\lambda;\uu_{-\infty})}\,,\\
\trho^\ell_\eps(\lambda)&\,=\,
D_\eps(\lambda)\,\frac{\eD^{-\int_{-\infty}^0(f'(\uU_\eps)-f'(\uu_{-\infty}))}}{\mu_+^\eps(\lambda;\uu_{+\infty})-\mu_-^\eps(\lambda;\uu_{+\infty})}\,.
\end{align*}
As a consequence for such a $(\eps,\lambda)$ and any $x\in\R$
\begin{align*}
\bV_\eps^{\ell,u}(\lambda,x)\cdot
\bp 0&-1\\1&0\ep
\tbV_\eps^{\ell,u}(\lambda,x)&=\,0\,,&
\bV_\eps^{r,s}(\lambda,x)\cdot
\bp 0&-1\\1&0\ep
\tbV_\eps^{r,s}(\lambda,x)&=\,0\,,
\end{align*}
\begin{align*}
\bV_\eps^{r,s}(\lambda,x)\cdot
\bp 0&-1\\1&0\ep
\tbV_\eps^{\ell,u}(\lambda,x)&=\,
-D_\eps(\lambda)\,\eD^{-\int_{-\infty}^0(f'(\uU_\eps)-f'(\uu_{-\infty}))}\,,\\
\bV_\eps^{\ell,u}(\lambda,x)\cdot
\bp 0&-1\\1&0\ep
\tbV_\eps^{r,s}(\lambda,x)&=\,
D_\eps(\lambda)\,\eD^{\int_0^{+\infty}(f'(\uU_\eps)-f'(\uu_{+\infty}))}\,.
\end{align*}
\epr

Notation $\tau$, $\rho$ is used here to echo transmission/reflection coefficients of the classical scattering framework.

A corresponding proposition holds for the high-frequency regime. 

\begin{proof}
All the properties are readily obtained by combining the fact that both $(\bV_\eps^{r,s}(\lambda,\cdot),\bV_\eps^{r,u}(\lambda,\cdot))$ and $(\bV_\eps^{\ell,s}(\lambda,\cdot),\bV_\eps^{\ell,u}(\lambda,\cdot))$ form a basis of solutions of the spectral system of ODEs, the Liouville formula for Wronskians and duality relation \eqref{eq:duality}.
\end{proof}
We thus have that : \begin{align*}
\bp \bV_\eps^{r,s}(y)&\bV_\eps^{\ell,u}(y)\ep\,^{-1}= \begin{pmatrix}
\dfrac{e_2 \cdot \tbV_\eps^{\ell,u}(\lambda,y)}{\,
D_\eps(\lambda)\,\eD^{-\int_{-\infty}^0(f'(\uU_\eps)-f'(\uu_{-\infty}))}\,} & \dfrac{e_1 \cdot \tbV_\eps^{\ell,u}(\lambda,y)}{\,
-D_\eps(\lambda)\,\eD^{-\int_{-\infty}^0(f'(\uU_\eps)-f'(\uu_{-\infty}))}\,} \\
-\dfrac{e_2\cdot \tbV_\eps^{r,s}(\lambda,y)}{\,
D_\eps(\lambda)\,\eD^{\int_0^{+\infty}(f'(\uU_\eps)-f'(\uu_{+\infty}))}\,} & \dfrac{e_1 \cdot \tbV_\eps^{r,s}(\lambda,y)}{\,
D_\eps(\lambda)\,\eD^{\int_0^{+\infty}(f'(\uU_\eps)-f'(\uu_{+\infty}))}\,}
\end{pmatrix}\,.
\end{align*}
\bpr\label{p:duality-high}
There exist positive constants $(\eps_0,C,\delta)$ such that setting with $\Omega_\delta$ as in Proposition~\ref{p:tracking-lemma} 
\[
\Omega_\delta:=\left\{\,\lambda\,;\ 
\Re\left(\sqrt{\lambda}\right)
\,\geq \frac{1}{\delta}
\right\}
\]
there exist on $[0,\eps_0]\times \Omega_\delta$ maps $\tau^{r,HF}$, $\rho^{r,HF}$, $\tau^{\ell,HF}$, $\rho^{\ell,HF}$, $\ttau^{r,HF}$, $\trho^{r,HF}$, $\ttau^{\ell,HF}$, $\trho^{\ell,HF}$, satisfying high-frequency versions of the conclusions of Proposition~\ref{p:duality-mean} and moreover all these functions are uniformly bounded on $[0,\eps_0]\times \Omega_\delta$.
\epr

\begin{proof}
Most of the proof is contained in the proof of Proposition~\ref{p:duality-mean}. The remaining part is directly derived from the observation that Proposition~\ref{p:tracking-lemma} provides asymptotics for $(P_\eps^\ell(\lambda,0))^{-1}\,P_\eps^r(\lambda,0)$ thus also for the coefficients under consideration.
\end{proof}

This leads to the following definition (wherever it makes sense)
\begin{align}\label{eq:spectral-Green}
G_\eps(\lambda;x,y)
\,:=\,
\begin{cases}
-\dfrac{\eD^{\int_{-\infty}^0(f'(\uU_\eps)-f'(\uu_{-\infty}))}}{D_\eps(\lambda)}\,
\beD_1\cdot\bV_\eps^{r,s}(\lambda,x)\,
\beD_1\cdot\tbV_\eps^{\ell,u}(\lambda,y)
&\quad\textrm{if }x>y\,,\\[0.5em]
-\dfrac{\eD^{-\int_0^{+\infty}(f'(\uU_\eps)-f'(\uu_{+\infty}))}}{D_\eps(\lambda)}
\beD_1\cdot\bV_\eps^{\ell,u}(\lambda,x)\,
\beD_1\cdot\tbV_\eps^{r,s}(\lambda,y)
&\quad\textrm{if }x<y\,,
\end{cases}
\end{align}
where $\beD_1:=(1,0)$. Note that $(\lambda-\cL_\eps)G_\eps(\lambda;\cdot,y)=\delta_y$ and, for $\Re(\lambda)$ sufficiently large, $G_\eps(\lambda;x,y)$ is exponentially decaying as $\|x-y\|\to\infty$. To bound $G_\eps(\lambda;x,y)$, we shall refine the alternative $x<y$ \emph{vs.} $x>y$. For instance, when $x>y$, more convenient equivalent representations of $G_\eps(\lambda;x,y)$ are
\begin{align*}
-\dfrac{\eD^{\int_{-\infty}^0(f'(\uU_\eps)-f'(\uu_{-\infty}))}}{D_\eps(\lambda)}\,
\beD_1\cdot\bV_\eps^{r,s}(\lambda,x)\,
\beD_1\cdot\tbV_\eps^{\ell,u}(\lambda,y)
&\quad\textrm{when }x>0>y\,,\\[0.5em]
-\eD^{\int_{-\infty}^0(f'(\uU_\eps)-f'(\uu_{-\infty}))}\,
\beD_1\cdot\bV_\eps^{r,s}(\lambda,x)\,
\beD_1\cdot\left(\dfrac{\trho^\ell_\eps(\lambda)}{D_\eps(\lambda)}\,\tbV_\eps^{r,u}(\lambda,y)
+\dfrac{\ttau^\ell_\eps(\lambda)}{D_\eps(\lambda)}\,\tbV_\eps^{r,s}(\lambda,y)\right)
&\quad\textrm{when }x>y>0\,,\\[0.5em]
-\eD^{\int_{-\infty}^0(f'(\uU_\eps)-f'(\uu_{-\infty}))}\,
\beD_1\cdot\left(
\dfrac{\rho^r_\eps(\lambda)}{D_\eps(\lambda)}\,\bV_\eps^{\ell,s}(\lambda,x)
+\dfrac{\tau^r_\eps(\lambda)}{D_\eps(\lambda)}\,\bV_\eps^{\ell,u}(\lambda,x)\right)\,
\beD_1\cdot\tbV_\eps^{\ell,u}(\lambda,y)
&\quad\textrm{when }0>x>y\,.
\end{align*}

\br
The representation of spectral Green functions, thus of resolvent operators, with Evans' functions is sufficient to prove classical results about the identification of spectrum --- including algebraic multiplicity --- at the right-hand side of the essential spectrum with zeros of Evans' functions.
\er

We use similar formulas in the high-frequency regime. Yet the Green functions of the high-frequency regime and the compact-frequency regime agree where they co-exist (by uniqueness of the spectral problem (in a suitably weighted space) in some overlapping regions and uniqueness of analytic continuation elsewhere) so that we do not need to introduce a specific piece of notation for the high-frequency regime.

We also point out that it follows from Proposition~\ref{p:tracking-lemma} that in the zone of interest 
\[
|\d_xG_\eps(\lambda;x,y)|\leq C\,\max(\{1,\sqrt{|\lambda|}\})\,|G_\eps(\lambda;x,y)|
\]
for some uniform constant $C$.

\subsection{Time-evolution}

It follows from standard semigroup theory that the representation 
\be\label{eq:semi-group}
S_\eps(t)\,=\,\frac{1}{2\iD\pi}\,\int_\R \eD^{\Lambda(\xi)\,t}\Lambda'(\xi)\,
(\Lambda(\xi)\I-\cL_\eps)^{-1}\,\dD\xi 
\ee
holds in $\cL(BUC^0(\R))$ when $\Lambda:\R\to \C$ is a continuous, piecewise $\cC^1$ simple curve such that
\begin{enumerate}
\item $\Lambda$ is valued in the right-hand connected component of\footnote{This set contains $\{\lambda;\Re(\lambda)\geq \omega\}$ when $\omega$ is sufficiently large.}
\[
\left\{\,\lambda\,;\ \textrm{for }u\in\{\uu_{-\infty},\uu_{+\infty}\}\,,\ 
\Re\left(\mu_+^\eps(\lambda,u)\right)>0>\Re\left(\mu_-^\eps(\lambda,u)\right)
\right\}\,;
\] 
\item there hold
\begin{align*}
\lim_{\xi\to\pm\infty}\Im(\Lambda(\xi))&\,=\,\pm\infty\,,&
\int_\R \eD^{\Re(\Lambda(\xi))\,t}
\frac{|\Lambda'(\xi)|}{1+|\Lambda(\xi)|}\dD\xi <+\infty\,
\end{align*}
and there exist positive $(R,c)$ such that for $|\xi|\geq R$
\[
\Re(\Lambda(\xi))\geq -c\,|\Im(\Lambda(\xi))|\,;
\]
\item there is no root of $D_\eps$ on the right\footnote{The second condition implies that this makes sense.} of $\Lambda(\R)$. 
\end{enumerate}
Failure of the third condition could be restored by adding positively-oriented small circles to the contour $\Lambda$. This is the first condition that we want to relax by going to Green functions.

For curves as above, applying the above formula to functions in $W^{\infty,\infty}(\R)$ and testing it against functions in $\cC^\infty_c(\R)$ leads to a similar representation for Green functions
\be\label{eq:time-Green}
G^\eps_t(x,y)\,=\,\frac{1}{2\iD\pi}\,\int_\R \eD^{\Lambda(\xi)\,t}\Lambda'(\xi)\,
G_\eps(\Lambda(\xi);x,y)\,\dD\xi\,.
\ee
The point is that at fixed $(t,x,y)$, the constraints on $\Lambda$ ensuring the representation formula are significantly less stringent and one may use this freedom to optimize bounds. In particular depending on the specific regime for the triplet $(t,x,y)$ or the kind of data one has in mind, one may trade spatial localization for time-decay and \emph{vice versa} by adjusting contours to the right so as to gain spatial decay or to the left in order to improve time decay. 

When doing so, we essentially follow the strategy of \cite{Zumbrun-Howard}. The critical decay is essentially encoded in limiting-endstates spectral spatial decay and Evans' function root location. Therefore, roughly speaking, leaving aside questions related to the presence of a root of the Evans' function at zero, contours are chosen here to approximately\footnote{In some cases a genuine optimization --- as in direct applications of the Riemann saddle point method --- would be impractical.} optimize bounds on
\[
\int_\R \eD^{\Re(\Lambda(\xi))\,t+\Re(\mu_\sharp^\eps(\Lambda(\xi),\uu_{\sign(x)\infty}))\,x
+\Re(\tmu_\flat^\eps(\Lambda(\xi),\uu_{\sign(y)\infty}))\,y}\frac{|\Lambda'(\xi)|}{|D_\eps(\Lambda(\xi))|}\,\dD\xi
\]
with $(\sharp,\flat)\in\{+,-\}^2$. More precisely, at fixed $(t,x,y)$, one picks $\Lambda_0$ real in $[-\tfrac12\eta_0,+\infty]$ (with $\eta_0$ as in Proposition~\ref{p:Sturm-Liouville}), approximately minimizing
\[
\Re(\lambda)\,t+\Re(\mu_\sharp^\eps(\lambda,\uu_{\sign(x)\infty}))\,x
+\Re(\tmu_\flat^\eps(\lambda,\uu_{\sign(y)\infty}))\,y
\]
among such real $\lambda$ in $[-\tfrac12\eta_0,+\infty]$ and then depending on cases one defines $\Lambda$ through one of the equations  
\begin{align*}
\Re(\mu_\sharp^\eps(\Lambda(\xi),\uu_{\sign(x)\infty}))\,x
&+\Re(\tmu_\flat^\eps(\Lambda(\xi),\uu_{\sign(y)\infty}))\,y\\
&=\Re(\mu_\sharp^\eps(\Lambda_0,\uu_{\sign(x)\infty}))\,x
+\Re(\tmu_\flat^\eps(\Lambda_0,\uu_{\sign(y)\infty}))\,y
+\iD\,\xi\,\zeta_{\sign(\xi)}\,(\sharp\,x+\flat\,y)\\
\Re(\tmu_\flat^\eps(\Lambda(\xi),\uu_{\sign(y)\infty}))\,y
&=\Re(\tmu_\flat^\eps(\Lambda_0,\uu_{\sign(y)\infty}))\,y
+\iD\,\xi\,\zeta_{\sign(\xi)}\,\times(\flat\,y)\\
\Re(\mu_\sharp^\eps(\Lambda(\xi),\uu_{\sign(x)\infty}))\,x
&=\Re(\mu_\sharp^\eps(\Lambda_0,\uu_{\sign(x)\infty}))\,x
+\iD\,\xi\,\zeta_{\sign(\xi)}\,\times(\sharp\,x)\\
\end{align*}
with $\zeta_\pm$ conveniently chosen to ensure a condition analogous to the second condition of the semigroup representation and including
\begin{align*}
\lim_{|\xi|\to\infty}\Re(\sqrt{\Lambda(\xi)})&=+\infty\,,&
\lim_{\xi\to\pm\infty}\Im(\Lambda(\xi))&=\pm\infty\,.
\end{align*}
This should be thought as an approximate/simplified version of the saddlepoint method in the sense that $\Lambda_0=\Lambda(0)$ is an approximate maximizer of the exponential decay rate among real numbers, but a minimizer along the curve $\Lambda(\cdot)$. 

Computational details --- carried out in next section --- are cumbersome but the process is rather systematic.

\section{Linear stability}\label{s:linear}

We now make the most of our spectral preparation to derive linear stability estimates.

To motivate the analysis, let us anticipate that our achievement is the splitting of $(S_\eps(t))_{t\geq0}$ as 
\be\label{eq:linear-phase}
S_\eps(t)(w)(x)\,=\,\uU_\eps'(x)\,\Spsi(t)(w)\,+\,\Stau(t)(w)(x)\,.
\ee
for some $(\Spsi(t))_{t\geq0}$, $(\Stau(t))_{t\geq0}$ with $\Stau(0)=\Id$, so that the following proposition holds. 

\bpr\label{p:linear}
There exists $\eps_0>0$ such that 
\begin{enumerate}
\item there exists $C>0$ such that for any $t\geq0$, any $0\leq\eps\leq\eps_0$ and any $w\in BUC^0(\R)$
\begin{align*}
\|\Stau(t)(w)\|_{L^\infty(\R)}&+\min\left(\{1,\sqrt{t}\}\right)
\|\d_x\Stau(t)(w)\|_{L^\infty(\R)}+|\d_t\Spsi(t)(w)|\\
&\leq\,C\,\eD^{-\min\left(\left\{|g'(\uu_{-\infty})|,|g'(\uu_{+\infty})|\right\}\right)\,\eps\,t}\,\|w\|_{L^\infty(\R)}\,,
\end{align*}
and when moreover $w\in BUC^1(\R)$
\begin{align*}
\|\Stau(t)(w)\|_{W^{1,\infty}(\R)}
&\leq\,C\,\eD^{-\min\left(\left\{|g'(\uu_{-\infty})|,|g'(\uu_{+\infty})|\right\}\right)\,\eps\,t}\,\|w\|_{W^{1,\infty}(\R)}\,,
\end{align*}
\item for any $\theta>0$ there exist positive $(C_\theta,\omega_\theta)$ such that for any $t\geq0$, any $0\leq\eps\leq\eps_0$ and any $w\in BUC^0(\R)$
\begin{align*}
\|\Stau(t)(w)\|_{L^\infty(\R)}&+\min\left(\{1,\sqrt{t}\}\right)
\|\d_x\Stau(t)(w)\|_{L^\infty(\R)}+|\d_t\Spsi(t)(w)|\\
&\leq\,C_\theta\,\eD^{-\omega_\theta\,t}\,\|\eD^{\theta\,|\,\cdot\,|}\,w\|_{L^\infty(\R)}\,.
\end{align*}
\end{enumerate}
\epr 

Estimates on operators are derived through pointwise bounds on Green kernels from the trivial fact that if $\bT$ is defined through
\[
\bT(w)(x)\,=\,
\int_\R \bK(x,y)\,w(y)\,\dD y
\]
then 
\[
\|\bT(w)\|_{L^\infty(\R)}
\,\leq\,\|\bK\|_{L^\infty_x(L^1_y)}\,\|w\|_{L^\infty(\R)}\,.
\]

\subsection{Auxiliary lemmas}

To begin with, to gain a practical grasp on the way the placement of spectral curves impacts decay rates, we provide two lemmas, that will be of ubiquitous use when establishing pointwise bounds on Green functions. 

Both lemmas are motivated by the fact that when $\beta\geq0$ and $t>0$ the minimization of
\[
\Lambda_0\,t+\left(\frac{\alpha}{2}-\,\sqrt{\frac{\alpha^2}{4}-\,b\,+\Lambda_0}\right)\,\beta
\]
over $\Lambda_0\in (-\frac{\alpha^2}{4}+\,b\,,+\infty)$ is equivalent to
\be\label{eq:natural}
2\,\sqrt{\frac{\alpha^2}{4}-\,b\,+\Lambda_0}
\,=\,\frac{\beta}{t}\,.
\ee
The first lemma directly elucidates the consequences of this choice of $\Lambda_0$ in the approximate saddlepoint method sketched above.

\bl\label{l:curves-large}
Let $t>0$, $\alpha\in\R$, $\beta\geq0$, $\beta_0\geq0$ $b<0$, and $(\zeta_-,\zeta_+)\in\C^2$ such that 
\begin{align*}
\Re(\zeta_\pm)&>|\Im(\zeta_\pm)|\,,&
\mp\Im(\zeta_\pm)>0\,.
\end{align*}
Then the curve $\Lambda\,:\ \R\to\C$ defined through\footnote{Sign conditions on $\Im(\zeta_\pm)$ ensure that this is a licit definition.}
\[
2\,\sqrt{\frac{\alpha^2}{4}-\,b\,+\Lambda(\xi)}
\,=\,\frac{\beta_0}{t}+\iD\xi\,\zeta_{\sign(\xi)}\,,
\]
satisfies for any $\xi\in\R$, when either $\beta=\beta_0$ or ($\beta\geq\beta_0$ and $\alpha\leq0$)
\begin{align*}
\Re\left(\Lambda(\xi)\,t+\left(\frac{\alpha}{2}-\,\sqrt{\frac{\alpha^2}{4}-\,b\,+\Lambda(\xi)}\right)\,\beta\right)
&\leq-\left(\frac{\alpha^2}{4}+\,|b|\right)\,t
+\frac{\alpha}{2}\,\beta_0-\frac{\beta_0^2}{4\,t}-\frac{\xi^2}{4}\,\Re(\zeta_{\sign(\xi)}^2)\,t\\
&\qquad=-\,|b|\,t
-\frac{(\beta_0-\alpha\,t)^2}{4\,t}-\frac{\xi^2}{4}\,\Re(\zeta_{\sign(\xi)}^2)\,t
\end{align*}
and for any $\xi\in\R^*$
\[
|\Lambda'(\xi)|\,\leq\,|\zeta_{\sign(\xi)}|\,\left(1+\frac{\Re(\zeta_{\sign(\xi)})}{|\Im(\zeta_{\sign(\xi)})|}\right)
\Re\left(\sqrt{\frac{\alpha^2}{4}-\,b\,+\Lambda(\xi)}\right)\,.
\]
\el

We omit the proof of Lemma~\ref{l:curves-large} as straightforward and elementary.

The second lemma is designed to deal with cases when the natural choice \eqref{eq:natural} is not available because of extra constraints arising from Evans' function possible annulation in $(-\frac{\alpha^2}{4}+\,b\,,0)$. Explicitly, it focuses on the case when $\beta/t\leq \omega_0$ when $\omega_0$ is typically picked as either $\omega_r^{\eta_0}$ or $\omega_r^{\eta_0}$ with
\begin{align}\label{eq:om0}
\omega_r^{\eta_0}&:=
2\,\sqrt{\frac{(f'(\uu_{+\infty})-\sigma_\eps)^2}{4}\,-\frac{\eta_0}{2}}\,,&
\omega_\ell^{\eta_0}&:=
2\,\sqrt{\frac{(f'(\uu_{-\infty})-\sigma_\eps)^2}{4}\,-\frac{\eta_0}{2}}\,,
\end{align}
where $\eta_0$ is as in Proposition~\ref{p:Sturm-Liouville}. Since $\beta/t\leq \omega_0$ should be thought as a bounded-domain restriction, it is useful to let the second lemma also encode the possible trade-off between spatial localization and time decay.

\bl\label{l:curves-small}
Let $t>0$, $\alpha\in\R$, $\beta\geq0$, $b<0$, $(\zeta_-,\zeta_+)\in\C^2$ such that 
\begin{align*}
\Re(\zeta_\pm)&>|\Im(\zeta_\pm)|\,,&
\mp\Im(\zeta_\pm)>0\,,
\end{align*}
and $\omega_0\geq 0$ such that
\[
\beta\,\leq\,\omega_0\,t\,.
\]
Then the curve $\Lambda\,:\ \R\to\C$ defined through
\[
2\,\sqrt{\frac{\alpha^2}{4}-\,b\,+\Lambda(\xi)}
\,=\,\omega_0+\iD\xi\,\zeta_{\sign(\xi)}\,,
\]
satisfies for any $\xi\in\R^*$
\[
|\Lambda'(\xi)|\,\leq\,|\zeta_{\sign(\xi)}|\,\left(1+\frac{\Re(\zeta_{\sign(\xi)})}{|\Im(\zeta_{\sign(\xi)})|}\right)
\Re\left(\sqrt{\frac{\alpha^2}{4}-\,b\,+\Lambda(\xi)}\right)\,.
\]
and for any $\xi\in\R$ and $\eta>0$,
\begin{align*}
&\Re\left(\Lambda(\xi)\,t+\left(\frac{\alpha}{2}-\,\sqrt{\frac{\alpha^2}{4}-\,b\,+\Lambda(\xi)}\right)\,\beta\right)\\
&\leq-\,\left(\frac{\alpha^2-(1+\eta)\,\omega_0^2}{4}+|b|\right)\,t
-\frac{\omega_0-\alpha}{2}\,\beta
-\frac{\xi^2}{4}\,\left(\Re(\zeta_{\sign(\xi)}^2)-\frac{1}{\eta}\,|\Im(\zeta_{\sign(\xi)})|^2\right)\,t
\end{align*}
and, when moreover $\omega_0<|\alpha|$,
\begin{align*}
&\Re\left(\Lambda(\xi)\,t+\left(\frac{\alpha}{2}-\,\sqrt{\frac{\alpha^2}{4}-\,b\,+\Lambda(\xi)}\right)\,\beta\right)\\
&\leq-\,|b|\,t
-\frac{(|\alpha|\,t-\beta)^2}{4\,t}\,\left(1-(1+\eta)\,\frac{\omega_0^2}{\alpha^2}\right)
-\frac{\xi^2}{4}\,\left(\Re(\zeta_{\sign(\xi)}^2)-\frac{1}{\eta}\,|\Im(\zeta_{\sign(\xi)})|^2\right)\,t\,.
\end{align*}
\el

Note that to guarantee for some $\eta>0$ both 
\begin{align*}
\Re(\zeta_{\sign(\xi)}^2)&>\frac{1}{\eta}\,|\Im(\zeta_{\sign(\xi)})|^2\,,&
1&>(1+\eta)\,\frac{\omega_0^2}{\alpha^2}\,,
\end{align*}
one needs to enforce
\begin{align}\label{eq:curve-condition}
\Re(\zeta_{\sign(\xi)})>\frac{|\alpha|}{\sqrt{\alpha^2-\omega_0^2}}\,|\Im(\zeta_{\sign(\xi)})|\,.
\end{align}
Likewise when $\omega_0>|\alpha|$, one may extract large-time decay for $|\xi|\geq \xi_0>0$ provided that $\Re(\zeta_{\sign(\xi)})$ is sufficiently large.

\begin{proof} The starting point is that for any $\eta>0$,
\begin{align*}
&\Re\left(\Lambda(\xi)\,t+\left(\frac{\alpha}{2}-\,\sqrt{\frac{\alpha^2}{4}-\,b\,+\Lambda(\xi)}\right)\,\beta\right)\\
&\leq-\,\left(\frac{\alpha^2-\omega_0^2\,-\eta\,\left(\omega_0-\frac{\beta}{t}\right)^2}{4}+|b|\right)\,t
-\frac{\omega_0-\alpha}{2}\,\beta
-\frac{\xi^2}{4}\,\left(\Re(\zeta_{\sign(\xi)}^2)-\frac{1}{\eta}\,|\Im(\zeta_{\sign(\xi)})|^2\right)\,t\\
&\qquad=
-\,|b|\,t
-\frac{(\alpha\,t-\beta)^2}{4\,t}\,
-\frac{\xi^2}{4}\,\left(\Re(\zeta_{\sign(\xi)}^2)-\frac{1}{\eta}\,|\Im(\zeta_{\sign(\xi)})|^2\right)\,t
+\frac{(\omega_0\,t-\beta)^2}{4\,t}\,\left(1+\eta\right)\,.
\end{align*}
The first bound on the real part is then obtained by using the first formulation of the foregoing bound jointly with 
\[
\left(\omega_0-\frac{\beta}{t}\right)^2\leq \omega_0^2
\]
whereas the second bound, specialized to the case $|\alpha|>\omega_0$, stems from the second formulation and
\begin{align*}
(|\alpha|\,t-\beta)^2&\leq\ (\alpha\,t-\beta)^2\,,&
(\omega_0\,t-\beta)^2&\leq\ \frac{\omega_0^2}{\alpha^2}\,(|\alpha|\,t-\beta)^2\,.
\end{align*}
\end{proof}

\subsection{First separations}

We would like to split $G^\eps_t(x,y)$ into pieces corresponding to different behaviors. Yet we must take into account that our description of $G_\eps(\lambda;x,y)$ is different in high-frequency and compact regimes. To do so, we pick some curves and break them into pieces. 

Explicitly, motivated by \eqref{eq:curve-condition} with $\omega_0$ either $\omega_r^{\eta_0}$ or $\omega_r^{\eta_0}$ --- defined in \eqref{eq:om0} with $\eta_0$ as in Proposition~\ref{p:Sturm-Liouville} ---, we first choose $\zeta^{HF}_\pm$ such that
\begin{align*}
\Re(\zeta^{HF}_\pm)&\geq 2\,\frac{\max\left(\left\{|f'(\uu_{+\infty})-\sigma_\eps|,f'(\uu_{-\infty})-\sigma_\eps\right\}\right)}{\sqrt{2\,\eta_0}}\,|\Im(\zeta^{HF}_\pm)|\,,&
\mp\Im(\zeta^{HF}_\pm)>0\,.
\end{align*}
Then we define curves $\Lambda_\eps^r$, $\Lambda_\eps^\ell$ through
\begin{align*}
2\,\sqrt{\frac{(f'(\uu_{+\infty})-\sigma_\eps)^2}{4}-\,\eps\,g'(\uu_{+\infty})\,+\Lambda_\eps^r(\xi)}
&\,=\,\omega_r^{HF}+\iD\xi\,\zeta^{HF}_{\sign(\xi)}\,,\\
2\,\sqrt{\frac{(f'(\uu_{-\infty})-\sigma_\eps)^2}{4}-\,\eps\,g'(\uu_{-\infty})\,+\Lambda_\eps^\ell(\xi)}
&\,=\,\omega_\ell^{HF}+\iD\xi\,\zeta^{HF}_{\sign(\xi)}\,,
\end{align*}
where $\omega_r^{HF}$ and $\omega_\ell^{HF}$ are fixed such that
\begin{align*}
\omega_r^{HF}&>|f'(\uu_{+\infty})-\sigma_0|\,,&
\omega_\ell^{HF}&>|f'(\uu_{-\infty})-\sigma_0|\,.
\end{align*}
Note that this is sufficient to guarantee that both curves satisfy requirements ensuring \eqref{eq:semi-group} thus also \eqref{eq:time-Green}.

We shall do a particular treatment of the parts of the curves corresponding to $|\xi|\leq\xi^{HF}$ where we choose $\xi^{HF}$ as
\begin{align*}
\xi^{HF}
:= \frac{2\,\max\left(\left\{\omega_r^{HF}-\omega_r^{\eta_0},\omega_\ell^{HF}-\omega_\ell^{\eta_0}\right\}\right)}{\min\left(\left\{|\Im(\zeta^{HF}_+)|,|\Im(\zeta^{HF}_-)|\right\}\right)}\,.
\end{align*}
Once again the motivation for the definition of $\xi^{HF}$ stems from Lemma~\ref{l:curves-small}. Indeed the definition ensures that for $\omega\in\R$, a curve $\Lambda$, defined through 
\begin{align*}
2\,\sqrt{\frac{(f'(\uu_{+\infty})-\sigma_\eps)^2}{4}-\,\eps\,g'(\uu_{+\infty})\,+\Lambda(\xi)}
&\,=\,\omega+\iD\xi\,\zeta^{\omega}_{r,\sign(\xi)}\,,
\end{align*}
respectively through
\begin{align*}
2\,\sqrt{\frac{(f'(\uu_{-\infty})-\sigma_\eps)^2}{4}-\,\eps\,g'(\uu_{-\infty})\,+\Lambda(\xi)}
&\,=\,\omega+\iD\xi\,\zeta^{\omega}_{\ell,\sign(\xi)}\,,
\end{align*}
with 
\[
\zeta^{\omega}_{\sharp,\pm}
:=\Re(\zeta^{HF}_\pm)+\iD\,\left(\Im\left(\zeta^{HF}_\pm\right)
\mp\frac{\omega_\sharp^{HF}-\omega}{\xi^{HF}}\right)\,,
\quad\sharp\in\{r,\ell\}\,,
\]
satisfies 
\begin{align*}
\Lambda(\pm\xi^{HF})&=\Lambda_\eps^r(\pm\xi^{HF})\,,&
\textrm{respectively }\quad
\Lambda(\pm\xi^{HF})&=\Lambda_\eps^\ell(\pm\xi^{HF})\,,
\end{align*}
whereas, for $\sharp\in\{r,\ell\}$, $\omega\in[\omega_\sharp^{\eta_0},\omega_\sharp^{HF}]$,
\begin{align*}
\Re(\zeta^{\omega}_{\sharp,\pm})&=\Re(\zeta^{HF}_\pm)\,,&
\mp\Im(\zeta^{\omega}_{\sharp,\pm})&>0\,,&
|\Im(\zeta^{\omega}_{\sharp,\pm})|&\leq \frac32|\Im(\zeta^{HF}_\pm)|\,.
\end{align*}
In the following, for $\sharp\in\{r,\ell\}$, we use notation $\Lambda_\eps^{\sharp,LF}:=(\Lambda_\eps^\sharp)_{|[-\xi^{HF},\xi^{HF}]}$ and $\Lambda_\eps^{\sharp,HF}:=(\Lambda_\eps^\sharp)_{|\R\setminus[-\xi^{HF},\xi^{HF}]}$. 

To ensure that Lemma~\ref{l:curves-small} provides exponential time decay for the part of the evolution arising from $\Lambda_\eps^{\sharp,HF}$, we reinforce the constraint on $\Re(\zeta^{HF}_\pm)$ by adding
\[
\Re(\zeta^{HF}_\pm)\ \geq \ 
\sqrt{\Im(\zeta^{HF}_\pm)^2+2\,\frac{\max\left(\left\{(\omega_r^{HF})^2-(f'(\uu_{+\infty})-\sigma_\eps)^2,
(\omega_\ell^{HF})^2-(f'(\uu_{-\infty})-\sigma_\eps)^2\right\}\right)}{(\xi^{HF})^2}}\,.
\]

Anticipating our needs when analyzing small-$\lambda$ expansions, we point out that by lowering $\eta_0$ and $\omega_\ell^{HF}$, $\omega_r^{HF}$, we may enforce that for $\sharp\in\{r,\ell\}$, when $\omega=\omega_\sharp^{\eta_0}$, and $\Lambda$ is defined as above, there exists $\omega'>0$ and some $\delta >0$, such that for any $\xi\in [-\xi^{HF},\xi^{HF}]$ 
\begin{align*}
\Re(\Lambda(\xi))&\leq -\omega'\,,&
\Re\left(\sqrt{\frac{(f'(\uu_{+\infty})-\sigma_\eps)^2}{4}-\,\eps\,g'(\uu_{+\infty})\,+\Lambda(\xi)}\right)\leq \frac{(-f'(\uu_{+\infty})+\sigma_\eps)}{2}\,- \delta \, ,
\end{align*}
respectively
\begin{align*}
\Re(\Lambda(\xi))&\leq -\omega'\,,&
\Re\left(\sqrt{\frac{(f'(\uu_{-\infty})-\sigma_\eps)^2}{4}-\,\eps\,g'(\uu_{-\infty})\,+\Lambda(\xi)}\right)\leq \frac{(f'(\uu_{-\infty})-\sigma_\eps)}{2}\, -\delta \,.
\end{align*}

After these preliminaries, to account for different behaviors, when $t>0$ we break $G^\eps_t$ as
\be
G^\eps_t\,=\,\Gtau_t+\Grho_t
\ee
with $\Gtau_t$ and $\Grho_t$ defined as follows.
First
\begin{align*}
\Gtau_t(x,y)
&=0\,,&\,\textrm{ if }xy>0\textrm{ and }
\left(\,y\geq\omega_r^{HF}\,t\ \textrm{or}\ 
y\leq-\omega_\ell^{HF}\,t\,\right)
\end{align*}
\begin{align*}
\Gtau_t(x,y)&\,=\,
\frac{1}{2\iD\pi}\,\int_\Lambda \eD^{\lambda\,t}\,
G_\eps(\lambda;x,y)\,\dD\lambda
&\,\textrm{if }xy<0\,,
\end{align*}
and when $xy>0$ and $-\omega_\ell^{HF}\,t<y<\omega_r^{HF}\,t$
\begin{align*}
&\Gtau_t(x,y)\\
&=
\begin{cases}
-\frac{1}{2\iD\pi}\,\int_{\Lambda_\eps^{r,LF}} \eD^{\lambda\,t}\,
\eD^{\int_{-\infty}^0(f'(\uU_\eps)-f'(\uu_{-\infty}))}\,
\beD_1\cdot\bV_\eps^{r,s}(\lambda,x)\,
\dfrac{\ttau^\ell_\eps(\lambda)}{D_\eps(\lambda)}\,
\beD_1\cdot\tbV_\eps^{r,s}(\lambda,y)\,\dD\lambda
&\,\textrm{if }x>y>0\\[0.5em]
-\frac{1}{2\iD\pi}\,\int_{\Lambda_\eps^{\ell,LF}} \eD^{\lambda\,t}\,
\eD^{\int_{-\infty}^0(f'(\uU_\eps)-f'(\uu_{-\infty}))}\,
\dfrac{\tau^r_\eps(\lambda)}{D_\eps(\lambda)}\,\beD_1\cdot\bV_\eps^{\ell,u}(\lambda,x)\,
\beD_1\cdot\tbV_\eps^{\ell,u}(\lambda,y)\,\dD\lambda
&\,\textrm{if }0>x>y\\[0.5em]
-\frac{1}{2\iD\pi}\,\int_{\Lambda_\eps^{\ell,LF}} \eD^{\lambda\,t}\,
\eD^{-\int_0^{+\infty}(f'(\uU_\eps)-f'(\uu_{+\infty}))}\,
\beD_1\cdot\bV_\eps^{\ell,u}(\lambda,x)\,
\dfrac{\ttau^r_\eps(\lambda)}{D_\eps(\lambda)}\,
\beD_1\cdot\tbV_\eps^{\ell,u}(\lambda,y)\,\dD\lambda
&\,\textrm{if }0>y>x\\[0.5em]
-\frac{1}{2\iD\pi}\,\int_{\Lambda_\eps^{r,LF}} \eD^{\lambda\,t}\,
\eD^{-\int_0^{+\infty}(f'(\uU_\eps)-f'(\uu_{+\infty}))}\,
\dfrac{\tau^\ell_\eps(\lambda)}{D_\eps(\lambda)}\,\beD_1\cdot\bV_\eps^{r,s}(\lambda,x)\,
\beD_1\cdot\tbV_\eps^{r,s}(\lambda,y)\,\dD\lambda
&\,\textrm{if }y>x>0\,,
\end{cases}
\end{align*}
where, here and in the definition of $\Grho_t$, $\Lambda$ is either $\Lambda=\Lambda_\eps^r$ or $\Lambda=\Lambda_\eps^\ell$, and we use compact notation for integrals over curves instead of explicitly parametrized versions. Second, 
\begin{align*}
\Grho_t(x,y)
&=0\,,&\,\textrm{ if }xy<0
\end{align*}
\begin{align*}
\Grho_t(x,y)
&=\frac{1}{2\iD\pi}\,\int_\Lambda \eD^{\lambda\,t}\,
G_\eps(\lambda;x,y)\,\dD\lambda\,,
&\,\textrm{ if }xy>0\textrm{ and }
\left(\,y\geq\omega_r^{HF}\,t\ \textrm{or}\ 
y\leq-\omega_\ell^{HF}\,t\,\right)\,,
\end{align*}
and when $xy>0$ and $-\omega_\ell^{HF}\,t<y<\omega_r^{HF}\,t$
\begin{align*}
&\Grho_t(x,y)\\
&=
\begin{cases}
\,-\frac{1}{2\iD\pi}\,\int_{\Lambda_\eps^{r,LF}} \eD^{\lambda\,t}\,
\eD^{\int_{-\infty}^0(f'(\uU_\eps)-f'(\uu_{-\infty}))}\,
\beD_1\cdot\bV_\eps^{r,s}(\lambda,x)\,
\dfrac{\trho^\ell_\eps(\lambda)}{D_\eps(\lambda)}\,
\beD_
1\cdot\tbV_\eps^{r,u}(\lambda,y)\,\dD\lambda&\\
\quad+\frac{1}{2\iD\pi}\,\int_{\Lambda_\eps^{r,HF}} \eD^{\lambda\,t}\,
G_\eps(\lambda;x,y)\,\dD\lambda
&\,\textrm{if }x>y>0\\[0.5em]
-\frac{1}{2\iD\pi}\,\int_{\Lambda_\eps^{\ell,LF}} \eD^{\lambda\,t}\,
\eD^{\int_{-\infty}^0(f'(\uU_\eps)-f'(\uu_{-\infty}))}\,
\dfrac{\rho^r_\eps(\lambda)}{D_\eps(\lambda)}\,\beD_1\cdot\bV_\eps^{\ell,s}(\lambda,x)\,
\beD_1\cdot\tbV_\eps^{\ell,u}(\lambda,y)\,\dD\lambda&\\
\quad+\frac{1}{2\iD\pi}\,\int_{\Lambda_\eps^{\ell,HF}} \eD^{\lambda\,t}\,
G_\eps(\lambda;x,y)\,\dD\lambda
&\,\textrm{if }0>x>y\\[0.5em]
-\frac{1}{2\iD\pi}\,\int_{\Lambda_\eps^{\ell,LF}} \eD^{\lambda\,t}\,
\eD^{-\int_0^{+\infty}(f'(\uU_\eps)-f'(\uu_{+\infty}))}\,
\beD_1\cdot\bV_\eps^{\ell,u}(\lambda,x)\,
\dfrac{\trho^r_\eps(\lambda)}{D_\eps(\lambda)}\,
\beD_1\cdot\tbV_\eps^{\ell,s}(\lambda,y)\,\dD\lambda&\\
\quad+\frac{1}{2\iD\pi}\,\int_{\Lambda_\eps^{\ell,HF}} \eD^{\lambda\,t}\,
G_\eps(\lambda;x,y)\,\dD\lambda
&\,\textrm{if }0>y>x\\[0.5em]
-\frac{1}{2\iD\pi}\,\int_{\Lambda_\eps^{r,LF}} \eD^{\lambda\,t}\,
\eD^{-\int_0^{+\infty}(f'(\uU_\eps)-f'(\uu_{+\infty}))}\,
\dfrac{\rho^\ell_\eps(\lambda)}{D_\eps(\lambda)}\,\beD_1\cdot\bV_\eps^{r,u}(\lambda,x)\,
\beD_1\cdot\tbV_\eps^{r,s}(\lambda,y)\,\dD\lambda&\\
\quad+\frac{1}{2\iD\pi}\,\int_{\Lambda_\eps^{r,HF}} \eD^{\lambda\,t}\,
G_\eps(\lambda;x,y)\,\dD\lambda
&\,\textrm{if }y>x>0\,.
\end{cases}
\end{align*} 
Note that the above splitting implies 
\[
\d_xG^\eps_t\,=\,\d_x\Gtau_t+\d_x\Grho_t
\]
where here and elsewhere throughout the text, $\d_x$ acting on either $\Gtau_t$ or $\Grho_t$ is understood as a pointwise derivative wherever these functions are continuous. 

The rationale behind the splitting is that the large-time decay of $\Grho_t$ is essentially limited by spatial decay hence may be thought as purely explained by essential spectrum considerations whereas the large-time asymptotics of $\Gtau_t$ is driven by the presence near the spectral curves of a root of $D_\eps$ at $\lambda=0$, hence is due to the interaction of essential and point spectra. 

Some extra complications in the splitting are due to the fact that we need to prepare the identification of the most singular part as a phase modulation, which comes into a tensorized form. This explains why we define zones in terms of the size of $|y|$, instead of the otherwise more natural $|x-y|$.

\subsection{First pointwise bounds}

We begin our use of Lemmas~\ref{l:curves-large} and~\ref{l:curves-small}  with short-time bounds.
\bl\label{l:short-time}
There exist positive $(\eps_0,C,\omega,\theta)$ such that for any $t>0$, any $0\leq\eps\leq\eps_0$ and any $(x,y)\in\R^2$
\begin{align*}
|\Gtau_t(x,y)|+\min\left(\{1,\sqrt{t}\}\right)|\d_x\Gtau_t(x,y)|&\,\leq\,C\,\eD^{\omega\,t}\,\eD^{-\theta\,|x|}\,\frac{1}{\sqrt{t}}\eD^{-\theta\,\frac{y^2}{t}}\,,
\end{align*}
\el

The foregoing lemma does not contain estimates on $\Grho_t$ because those would be redundant with the corresponding large-time estimates. The point of Lemma~\ref{l:short-time} is to show that for short-time estimates the singularity at $\lambda=0$ may be avoided whereas this singularity is not present in $\Grho_t$.

\begin{proof}
To bound $\Gtau_t(x,y)$ when $xy<0$, we separate between $x>0>y$ and $x<0<y$. The analyses being completely similar, we only discuss here the former case. To treat it, we move curves as in Lemmas~\ref{l:curves-large} and~\ref{l:curves-small} with $\beta_0=\beta=|y|$, $\alpha=f'(\uu_{-\infty})-\sigma_\eps$, $b=\eps\,g'(\uu_{-\infty})$ and note that
\[
\Re(\mu_-^\eps(\lambda,\uu_{+\infty}))\,\leq\,
\frac12\left(f'(\uu_{+\infty})-\sigma_\eps\right)
\,<\,0\,.
\] 
More explicitly, we use Lemma~\ref{l:curves-large} to bound the regime $|y|\geq \omega_\ell^{HF}\,t$ which leads to the claimed heat-like bound since
\begin{align*}
\eD^{-\frac{(|y|-\alpha\,t)^2}{4\,t}}
&\leq \eD^{-\left(1-\frac{|\alpha|}{\omega_\ell^{HF}}\right)\,\frac{y^2}{4\,t}}\,,&
&|y|\geq \omega_\ell^{HF}\,t\,.
\end{align*}
In the remaining zone where $|y|\leq \omega_\ell^{HF}\,t$ we use instead Lemma~\ref{l:curves-small} to derive a bound that may be converted into a heat-like bound through 
\begin{align*}
\eD^{-\frac{\omega_\ell^{HF}-\alpha}{2}|y|}
&\leq \eD^{-\left(1-\frac{|\alpha|}{\omega_\ell^{HF}}\right)\,\frac{y^2}{2\,t}}\,,&
&|y|\leq \omega_\ell^{HF}\,t\,.
\end{align*}

The estimates on $\Gtau_t(x,y)$ when $xy>0$ are obtained in exactly the same way.
\end{proof}


We proceed with bounds on $\Grho_t$.

\bl\label{l:Grho}
There exist positive $(\eps_0,C,\omega,\theta)$ such that for any $t>0$, any $0\leq\eps\leq\eps_0$ and any $(x,y)\in\R^2$
\begin{align*}
&|\Grho_t(x,y)|+\min\left(\{1,\sqrt{t}\}\right)|\d_x\Grho_t(x,y)|\\
&\quad\leq\,C\,\One_{|x-y|\leq |y|}\,\eD^{-\min\left(\left\{|g'(\uu_{-\infty})|,|g'(\uu_{+\infty})|\right\}\right)\,\eps\,t}\,\frac{1}{\sqrt{t}}\left(\eD^{-\theta\,\frac{|x-y-(f'(\uu_{+\infty})-\sigma_\eps)\,t|^2}{t}}
+\eD^{-\theta\,\frac{|x-y-(f'(\uu_{-\infty})-\sigma_\eps)\,t|^2}{t}}\right)\\
&\qquad+\,C\,\One_{|x-y|\geq |y|}\eD^{-\omega\,t}\,\frac{1}{\sqrt{t}}\eD^{-\theta\,\frac{y^2}{t}}\,.
\end{align*}
This also implies that there exist positive $(\eps_0,C,\theta)$ such that for any $\theta'>0$ there exists $\omega'>0$ such that for any $t>0$, any $0\leq\eps\leq\eps_0$ and any $(x,y)\in\R^2$
\begin{align*}
\eD^{-\theta'|y|}\,\Big(|\Grho_t(x,y)|&+\min\left(\{1,\sqrt{t}\}\right)|\d_x\Grho_t(x,y)|\Big)\\
&\quad\leq\,C\,\eD^{-\omega'\,t}\,\frac{1}{\sqrt{t}}\left(\eD^{-\theta\,\frac{|x-y-(f'(\uu_{+\infty})-\sigma_\eps)\,t|^2}{t}}
+\eD^{-\theta\,\frac{|x-y-(f'(\uu_{-\infty})-\sigma_\eps)\,t|^2}{t}}
+\eD^{-\theta\,\frac{y^2}{t}}\right)\,.
\end{align*}
\el

In the foregoing statement and throughout the text we use $\One_A$ to denote a characteristic function for the condition $A$.

\begin{proof}
To deduce the second bound from the first we observe that for any $\alpha$
\[
\frac{\theta}{2\,t}\,|x-y-\alpha\,t|^2+\theta'|x-y|
\geq\,\begin{cases}
\frac12\,\theta'\,|\alpha|\,t
&\,\textrm{if }|x-y-\alpha\,t|\leq \frac12|\alpha|\,t\\
\frac{\theta}{4}\,|\alpha|^2\,t
&\,\textrm{if }|x-y-\alpha\,t|\geq \frac12|\alpha|\,t
\end{cases}\,.
\]

To prove the first bound we should distinguish between regimes defined by $0<y<x$, $y<x<0$, $0>y>x$ and $y>x>0$. Regimes $0<y<x$ and $0>y>x$ on one hand and $y<x<0$ and $y>x>0$ on the other hand may be treated similarly and we give details only for the cases $y<x<0$ and $0<y<x$. 

Note that when $y<x<0$, we have $|x-y|\leq |y|$. When $y<x<0$ and $|y|\geq\omega_\ell^{HF}\,t$, we move the curve according to Lemma~\ref{l:curves-large} with $\beta_0=\beta=|x-y|$, $\alpha=f'(\uu_{-\infty})-\sigma_\eps$, $b=\eps\,g'(\uu_{-\infty})$. To analyze the regime when $y<x<0$ and $|y|<\omega_\ell^{HF}\,t$, we never move the curve $\Lambda_\eps^{\ell,HF}$ (but bound its contribution according to Lemma~\ref{l:curves-small}) whereas we move $\Lambda_\eps^{\ell,LF}$ as in Lemma~\ref{l:curves-large} with $\zeta_{\pm}=\zeta^{|x-y|/t}_{\ell,\pm}$ when $|x-y|\geq \omega_\ell^{\eta_0}\,t$ or as in Lemma~\ref{l:curves-small} with $\zeta_{\pm}=\zeta^{\omega_\ell^{\eta_0}}_{\ell,\pm}$ when $|x-y|\leq\omega_\ell^{\eta_0}\,t$.

To bound the contribution of the regime $0<y<x$, we may proceed as when $y<x<0$ provided that $|x-y|\leq |y|$ or $-\omega_\ell^{HF}\,t\leq y\leq \omega_r^{HF}\,t$. The remaining case is dealt with by applying Lemma~\ref{l:curves-large} with $\beta_0=|y|$ and $\beta=|x-y|$ using the fact that $f'(\uu_{+\infty})-\sigma_\eps<0$.
\end{proof}

\subsection{Linear phase separation}

The large-time estimates for $\Gtau_t$ require a phase separation. To carry it out we first recall that there exist $(a^r_\eps,a^\ell_\eps)\in\R^2$, each uniformly bounded from below and above, such that
\begin{align*}
\beD_1\cdot\bV_\eps^{r,s}(0,\cdot)&=a^r_\eps\,\uU'_\eps\,,&
\beD_1\cdot\bV_\eps^{\ell,u}(0,\cdot)&=a^\ell_\eps\uU'_\eps\,.
\end{align*}

Then we split $\Gtau_t$ as
\[
\Gtau_t(x,y)\,=\,\uU'_\eps(x)\,\Gpsi_t(y)\,+\,\tGtau_t(x,y)
\]
with 
\begin{align*}
\Gpsi_t(y)
&=0\,,&\,\textrm{ if }
\left(\,y\geq\omega_r^{HF}\,t\ \textrm{or}\ 
y\leq-\omega_\ell^{HF}\,t\,\right)\,,
\end{align*}
whereas when $-\omega_\ell^{HF}\,t<y<\omega_r^{HF}\,t$
\begin{align*}
\Gpsi_t(y)
&=
\begin{cases}\ds
-\frac{1}{2\iD\pi}\,\int_{\Lambda_\eps^{r,LF}} \eD^{\lambda\,t}\,
\eD^{-\int_0^{+\infty}(f'(\uU_\eps)-f'(\uu_{+\infty}))}\,
a^\ell_\eps\,
\beD_1\cdot\tbV_\eps^{r,s}(\lambda,y)\,\dfrac{\dD\lambda}{D_\eps(\lambda)}
&\,\textrm{if }y>0\\[1em]\ds
-\frac{1}{2\iD\pi}\,\int_{\Lambda_\eps^{\ell,LF}} \eD^{\lambda\,t}\,
\eD^{\int_{-\infty}^0(f'(\uU_\eps)-f'(\uu_{-\infty}))}\,
a^r_\eps\,
\beD_1\cdot\tbV_\eps^{\ell,u}(\lambda,y)
\dfrac{\dD\lambda}{D_\eps(\lambda)}
&\,\textrm{if }y<0\,.
\end{cases}
\end{align*}
As a result, when ($y\geq\omega_r^{HF}\,t$ or $y\leq-\omega_\ell^{HF}\,t$),
\begin{align*}
\tGtau_t(x,y)
&=
\begin{cases}
0\,,&\,\textrm{ if }xy>0\\
\frac{1}{2\iD\pi}\,\int_\Lambda \eD^{\lambda\,t}\,
G_\eps(\lambda;x,y)\,\dD\lambda
&\,\textrm{if }xy<0\,,
\end{cases}
\end{align*}
whereas when $xy<0$ and $-\omega_\ell^{HF}\,t<y<\omega_r^{HF}\,t$
\begin{align*}
&\tGtau_t(x,y)\\
&=
\begin{cases}
\,-\frac{1}{2\iD\pi}\,\int_{\Lambda_\eps^{r,LF}} \eD^{\lambda\,t}\,
\eD^{\int_{-\infty}^0(f'(\uU_\eps)-f'(\uu_{-\infty}))}\,
\beD_1\cdot\left(\bV_\eps^{r,s}(\lambda,x)-\bV_\eps^{r,s}(0,x)\right)\,
\beD_1\cdot\tbV_\eps^{\ell,u}(\lambda,y)\,\dfrac{\dD\lambda}{D_\eps(\lambda)}&\\
\quad+\frac{1}{2\iD\pi}\,\int_{\Lambda_\eps^{r,HF}} \eD^{\lambda\,t}\,
G_\eps(\lambda;x,y)\,\dD\lambda
\hfill\,\textrm{if }x>0>y\\[0.75em]
-\frac{1}{2\iD\pi}\,\int_{\Lambda_\eps^{\ell,LF}} \eD^{\lambda\,t}\,
\eD^{-\int_0^{+\infty}(f'(\uU_\eps)-f'(\uu_{+\infty}))}\,
\beD_1\cdot\left(\bV_\eps^{\ell,u}(\lambda,x)-\bV_\eps^{\ell,u}(0,x)\right)\,
\beD_1\cdot\tbV_\eps^{r,s}(\lambda,y)\,
\,\dfrac{\dD\lambda}{D_\eps(\lambda)}&\\
\quad+\frac{1}{2\iD\pi}\,\int_{\Lambda_\eps^{\ell,HF}} \eD^{\lambda\,t}\,
G_\eps(\lambda;x,y)\,\dD\lambda
\hfill\,\textrm{if }y>0>x\,,
\end{cases}
\end{align*} 
and when $xy>0$ and $-\omega_\ell^{HF}\,t<y<\omega_r^{HF}\,t$, $\tGtau_t(x,y)$ equals
\begin{align*}
\begin{cases}
-\frac{1}{2\iD\pi}\,\int_{\Lambda_\eps^{r,LF}} \eD^{\lambda\,t}\,
\eD^{\int_{-\infty}^0(f'(\uU_\eps)-f'(\uu_{-\infty}))}\,
\left(\ttau^\ell_\eps(\lambda)\beD_1\cdot\bV_\eps^{r,s}(\lambda,x)-\ttau^\ell_\eps(0)\beD_1\cdot\bV_\eps^{r,s}(0,x)\right)\,
\dfrac{1}{D_\eps(\lambda)}\,
\beD_1\cdot\tbV_\eps^{r,s}(\lambda,y)\,\dD\lambda
\\
\qquad-\frac{1}{2\iD\pi}\,\int_{\Lambda_\eps^{r,LF}} \eD^{\lambda\,t}\,
\eD^{\int_{-\infty}^0(f'(\uU_\eps)-f'(\uu_{-\infty}))}\,
\beD_1\cdot\bV_\eps^{r,s}(0,x)\,
\dfrac{\trho^\ell_\eps(\lambda)}{D_\eps(\lambda)}\,
\beD_1\cdot\tbV_\eps^{r,u}(\lambda,y)\,\dD\lambda
\hfill\,\textrm{if }x>y>0\\[0.5em]
-\frac{1}{2\iD\pi}\,\int_{\Lambda_\eps^{\ell,LF}} \eD^{\lambda\,t}\,
\eD^{\int_{-\infty}^0(f'(\uU_\eps)-f'(\uu_{-\infty}))}\,
\left(\tau^r_\eps(\lambda)\,\beD_1\cdot\bV_\eps^{\ell,u}(\lambda,x)
-\tau^r_\eps(0)\,\beD_1\cdot\bV_\eps^{\ell,u}(0,x)
\right)\,
\beD_1\cdot\tbV_\eps^{\ell,u}(\lambda,y)
\,\dfrac{\dD\lambda}{D_\eps(\lambda)}\\
\hfill\,\textrm{if }0>x>y\\[0.5em]
-\frac{1}{2\iD\pi}\,\int_{\Lambda_\eps^{\ell,LF}} \eD^{\lambda\,t}\,
\eD^{-\int_0^{+\infty}(f'(\uU_\eps)-f'(\uu_{+\infty}))}\,
\left(\ttau^r_\eps(\lambda)\,\beD_1\cdot\bV_\eps^{\ell,u}(\lambda,x)-\ttau^r_\eps(0)\,\beD_1\cdot\bV_\eps^{\ell,u}(0,x)\right)\,
\beD_1\cdot\tbV_\eps^{\ell,u}(\lambda,y)
\,\dfrac{\dD\lambda}{D_\eps(\lambda)}\\
\qquad-\frac{1}{2\iD\pi}\,\int_{\Lambda_\eps^{\ell,LF}} \eD^{\lambda\,t}\,
\eD^{-\int_0^{+\infty}(f'(\uU_\eps)-f'(\uu_{+\infty}))}\,
\beD_1\cdot\bV_\eps^{\ell,u}(0,x)\,
\dfrac{\trho^r_\eps(\lambda)}{D_\eps(\lambda)}\,
\beD_1\cdot\tbV_\eps^{\ell,s}(\lambda,y)\,\dD\lambda
\hfill\,\textrm{if }0>y>x\\[0.5em]
-\frac{1}{2\iD\pi}\,\int_{\Lambda_\eps^{r,LF}} \eD^{\lambda\,t}\,
\eD^{-\int_0^{+\infty}(f'(\uU_\eps)-f'(\uu_{+\infty}))}\,
\left(\tau^\ell_\eps(\lambda)\,\beD_1\cdot\bV_\eps^{r,s}(\lambda,x)
-\tau^\ell_\eps(0)\,\beD_1\cdot\bV_\eps^{r,s}(0,x)\right)\,
\beD_1\cdot\tbV_\eps^{r,s}(\lambda,y)\,\dfrac{\dD\lambda}{D_\eps(\lambda)}\\
\hfill\,\textrm{if }y>x>0\,,
\end{cases}
\end{align*}
where, here again $\Lambda$ is either $\Lambda=\Lambda_\eps^r$ or $\Lambda=\Lambda_\eps^\ell$.

In $\tGtau_t$ the contributions due to $\trho$-terms do not fit directly the framework of Lemmas~\ref{l:curves-large} and~\ref{l:curves-small}. This is here that instead we use that $\eta_0$, $\omega_r^{HF}$ and $\omega_\ell^{HF}$ were taken sufficient small to guarantee even simpler bounds when $\Lambda_\eps^{\sharp,LF}$ is moved as in Lemma~\ref{l:curves-small} with $\omega_0^\sharp=\omega_\sharp^{\eta_0}$, $\zeta_{\pm}=\zeta^{\omega_0^\sharp}_{\sharp,\pm}$, though we do not restrict to the zone $|y|\leq \omega_0^\sharp\,t$. Proceeding as above for the rest of bounds, we obtain the following lemmas.

\bl\label{l:tGtau}
There exist positive $(\eps_0,C,\omega,\theta)$ such that for any $t>0$, any $0\leq\eps\leq\eps_0$ and any $(x,y)\in\R^2$
\begin{align*}
&|\tGtau_t(x,y)|+\min\left(\{1,\sqrt{t}\}\right)|\d_x\tGtau_t(x,y)|\\
&\quad\leq\,C\,\eD^{-\min\left(\left\{|g'(\uu_{-\infty})|,|g'(\uu_{+\infty})|\right\}\right)\,\eps\,t}\,\eD^{-\theta\,|x|}\,
\frac{1}{\sqrt{t}}\left(\eD^{-\theta\,\frac{|y-(f'(\uu_{+\infty})-\sigma_\eps)\,t|^2}{t}}
+\eD^{-\theta\,\frac{|y+(f'(\uu_{-\infty})-\sigma_\eps)\,t|^2}{t}}\right)\\
&\qquad
+C\,\eD^{-\omega\,t}\,\eD^{-\theta\,|x|}\,\eD^{-\theta\,|x-y|}\,.
\end{align*}
This also implies that there exist positive $(\eps_0,C,\theta)$ such that for any $\theta'>0$ there exists $\omega'>0$ such that for any $t>0$, any $0\leq\eps\leq\eps_0$ and any $(x,y)\in\R^2$
\begin{align*}
\eD^{-\theta'|y|}\,&\Big(|\tGtau_t(x,y)|+\min\left(\{1,\sqrt{t}\}\right)|\d_x\tGtau_t(x,y)|\Big)\\
&\leq\,C\,\eD^{-\omega'\,t}\,\left(\frac{1}{\sqrt{t}}\left(\eD^{-\theta\,\frac{|y-(f'(\uu_{+\infty})-\sigma_\eps)\,t|^2}{t}}
+\eD^{-\theta\,\frac{|y+(f'(\uu_{-\infty})-\sigma_\eps)\,t|^2}{t}}\right)
+\,\eD^{-\theta\,|x|}\,\eD^{-\theta\,|x-y|}\,\eD^{-\theta'|y|}\right)\,.
\end{align*}
\el

\bl\label{l:Gpsi}
There exist positive $(\eps_0,C,\omega,\theta)$ such that for any $t>0$, any $0\leq\eps\leq\eps_0$ and any $y\in\R$
\begin{align*}
|\Gpsi_t(y)|&\,\leq\,C\,\eD^{\omega\,t}\,\frac{1}{\sqrt{t}}\eD^{-\theta\,\frac{y^2}{t}}\,.
\end{align*}
Moreover, there exist positive $(\eps_0,C,\omega,\theta)$ such that for any $t>0$, any $0\leq\eps\leq\eps_0$ and any $y\in\R$
\begin{align*}
&|\d_t\Gpsi_t(y)|
&\leq\,C\,\eD^{-\min\left(\left\{|g'(\uu_{-\infty})|,|g'(\uu_{+\infty})|\right\}\right)\,\eps\,t}\,
\frac{1}{\sqrt{t}}\left(\eD^{-\theta\,\frac{|y-(f'(\uu_{+\infty})-\sigma_\eps)\,t|^2}{t}}
+\eD^{-\theta\,\frac{|y+(f'(\uu_{-\infty})-\sigma_\eps)\,t|^2}{t}}\right)\,.
\end{align*}
This also implies that there exist positive $(\eps_0,C,\theta)$ such that for any $\theta'>0$ there exists $\omega'>0$ such that for any $t>0$, any $0\leq\eps\leq\eps_0$ and any $y\in\R$
\begin{align*}
\eD^{-\theta'|y|}\,|\d_t\Gpsi_t(y)|
&\leq\,C\,\eD^{-\omega'\,t}\,\frac{1}{\sqrt{t}}\left(\eD^{-\theta\,\frac{|y-(f'(\uu_{+\infty})-\sigma_\eps)\,t|^2}{t}}
+\eD^{-\theta\,\frac{|y+(f'(\uu_{-\infty})-\sigma_\eps)\,t|^2}{t}}\right)\,.
\end{align*}
\el

\hspace{2em}

To conclude and prove Proposition~\ref{p:linear}, we pick a smooth cut-off function $\chi$ on $[0,+\infty)$ such that $\chi\equiv 1$ on $[2,+\infty)$ and $\chi\equiv 0$ on $[0,1]$ and define $\Spsi$, $\Stau$ by
\begin{align*}
\Spsi(t)(w)&\,:=\,
\int_\R \chi(t)\,\Gpsi_t(y)\,w(y)\,\dD y\,,\\
\Stau(t)(w)(x)&\,:=\,
\int_\R \left(\chi(t)\,(\Grho_t(x,y)+\tGtau_t(x,y))+(1-\chi(t))\,G^\eps_t(x,y)\right)\,w(y)\,\dD y\,.
\end{align*}
The definitions are extended to $t=0$ by $\Spsi(0)(w)=0$ and $\Stau(0)(w)=w$.

As explained near its statement, Proposition~\ref{p:linear} follows then from $L^\infty_xL^1_y$ bounds on Green kernels, which themselves are derived from pointwise bounds proved above. 

\section{Nonlinear stability}

In the present section we conclude the proof of Theorem~\ref{th:main-sl}.

To do so, we seek for $u$ solving \eqref{eq:scalar-sl} under the form
\be\label{eq:stab-ansatz}
u(t,x)\,=\,\uU_\eps(x+\psi(t))+v(t,x+\psi(t))\,
\ee
with $(v,\psi')$ exponentially decaying in time. In these terms the equation becomes
\begin{align}\label{eq:scalar-ansatz}
\d_tv&+(f'(\uU_\eps+v)-\sigma_\eps+\psi')\d_x v-\d_x^2 v\\
&\nonumber
\,=\,\eps\,(g(\uU_\eps+v)-g(\uU_\eps))
-(f'(\uU_\eps+v)-f'(\uU_\eps)+\psi')\uU_\eps'\,.
\end{align}

Equation~\eqref{eq:scalar-ansatz} may be solved through
\begin{align}\label{eq:v}
v(t,\cdot)&\,=\,
\Stau(t)\,(v_0)\,+\,\int_0^t\Stau(t-s)\,\cN_\eps[v(s,\cdot)),\psi'(s)]\,\dD s\,,\\
\label{eq:psi}
\psi'(t)&\,=\,
\d_t\Spsi(t)\,(v_0)\,+\,\int_0^t\d_t\Spsi(t-s)\,\cN_\eps[v(s,\cdot)),\psi'(s)]\,\dD s\,,
\end{align}
with 
\begin{align*}
\cN_\eps[w,\varphi]&:=
-(f'(\uU_\eps+w)-f'(\uU_\eps)+\varphi)\d_x w
+\eps\,(g(\uU_\eps+w)-g(\uU_\eps)-g'(\uU_\eps)\,w)\\
&\quad-(f'(\uU_\eps+w)-f'(\uU_\eps)-f''(\uU_\eps)\,w)\uU_\eps'\,.
\end{align*}

In the present section, for notational concision's sake, we denote
\[
\omega_\infty=\min(|g'(\uu_{-\infty})|,|g'(\uu_{+\infty})|)\,.
\]

To begin with, we observe that estimates of the foregoing section are almost sufficient to run a continuity argument on \eqref{eq:v}-\eqref{eq:psi}. Indeed they provide the following proposition.

\bpr\label{p:Duhamel}
There exist $\theta_0>0$ and $\eps_0>0$ such that for any $0<\theta\leq \theta_0$ and $\delta>0$, there exist $C>0$ such that for any $0<\eps\leq\eps_0$ and $T>0$, if $(v,\psi')$ solves \eqref{eq:v}-\eqref{eq:psi} on $[0,T]$, with
\begin{align*}
\|v(t,\cdot)\|_{L^\infty(\R)}&\leq \delta\,,&\qquad t\in[0,T]\,,
\end{align*}
then, for any $t\in [0,T]$,
\begin{align*}
|\psi'(t)|&+\|v(t,\cdot)\|_{W^{1,\infty}(\R)}\\
&\leq\,
C\,\|v(0,\cdot)\|_{W^{1,\infty}(\R)}\,\eD^{-\eps\,\omega_\infty t}\\
&\quad\times\,\exp\left(
C\,\sup_{0\leq s\leq T}\eD^{\eps\,\omega_\infty s}(|\psi'(s)|
+\|v(s,\cdot)\|_{L^\infty(\R)}
+\|(\eps+\eD^{-\theta |\cdot|})^{-1}\d_xv(s,\cdot)\|_{L^\infty(\R)})\right)\,.
\end{align*}
\epr
 
The estimate fails to close by the fact that $\|\d_x w\|_{L^\infty(\R)}$ provides a weaker $\eps$-uniform control on $w$ than $\|(\eps+\eD^{-\theta |\cdot|})^{-1}\d_x w\|_{L^\infty(\R)}$. Note however that for any $x_*>0$, 
\[
\|(\eps+\eD^{-\theta |\cdot|})^{-1}\d_x w\|_{L^\infty([-x_*,x_*])}
\leq \eD^{\theta x_*}\,\|\d_x w\|_{L^\infty(\R)}
\]
so that we only need to improve the estimates on $\d_x v(t,\cdot)$ on the complement of some compact neighborhood of $0$.

\subsection{Maximum principle and propagation of regularity}\label{s:maximum}

To close our nonlinear estimates without using neither localization nor parabolic smoothing --- which would cause loss in powers of $\eps$---, we shall use a maximum principle argument. 

To begin with, we state and prove a convenient classical abstract maximum principle. We provide a proof mostly to highlight that it may be thought as an energy estimate on a suitable nonlinear function.

\bl\label{l:maximum}
Let $T>0$, $x_*\in\R$, $a\in L^1([0,T];W^{1,\infty}([x_*,+\infty)))$ bounded from above away from zero and $h\in\cC^0([0,T]\times[x_*,+\infty)\times\R)$. If $w\in\cC^2((0,T)\times[x_*,+\infty))\cap\cC^0([0,T]\times[x_*,+\infty))$ is a bounded function such that
\[
\d_t w+a(\cdot,\cdot)\,\d_x w\leq\d_x^2 w+h(\cdot,\cdot,w)\,,\qquad
\textrm{on }[0,T]\times[x_*,+\infty)
\]
and $M$ is a positive constant such that 
\begin{align*}
M&\geq w(\cdot,x_*)\,,&\qquad\textrm{on }[0,T]\,,\\
M&\geq w(0,\cdot,)\,,&\qquad\textrm{on }[x_*,+\infty)\,,\\
0&\geq \One_{\cdot>M}\,h(t,x,\cdot)\,,
\end{align*}
then 
\[
w\leq M\,,\qquad\textrm{on }[0,T]\times[x_*,+\infty)\,.
\]
\el

\begin{proof}
When moreover
\begin{align*}
M&> \limsup_{x\to\infty}w(\cdot,x)\,,&\qquad\textrm{on }[0,T]
\end{align*}
the claim is proved by a Gr\"onwall argument on
\[
t\mapsto\int_{x_*}^{+\infty}(w(t,x)-M)_+\,\dD x\,.
\]
The general case is recovered by applying this special case to $(t,x)\mapsto\eD^{-\theta\,(x-x_*)\,}w(t,x)$ with $\theta>0$ sufficiently small and taking the limit $\theta\to0$.
\end{proof}

We now use the foregoing lemma to derive a weighted bound on $\d_xv$ outside a sufficiently large compact neighborhood of $0$. We shall insert such a bound in a continuity argument so that we only need to prove that as long as $\d_xv$ does not become too large it remains small. This is the content of the following proposition.

\bpr\label{p:energy}
There exists $\theta_0>0$ such that for any $0<\theta\leq \theta_0$, there exist $x_*>0$, $\eps_0>0$, $\delta>0$ and $C>0$ such that for any $0<\eps\leq\eps_0$ and $T>0$, if $(v,\psi')$ solves \eqref{eq:scalar-ansatz} on $[0,T]\times\R$, with
\begin{align*}
|\psi'(t)|+\|v(t,\cdot)\|_{L^\infty(\R)}&\leq \delta \eD^{-\eps\,\omega_\infty t}\,,&\qquad t\in[0,T]\,,\\
\frac{|\d_xv(t,x)|}{\eps+\eD^{-\theta |x|}}&\leq \delta \eD^{-\eps\,\omega_\infty t}\,,&\qquad (t,x)\in [0,T]\times\R\,,
\end{align*}
then for any $(t,x)\in [0,T]\times(\R\setminus [-x_*,x_*])$
\begin{align*}
\frac{|\d_xv(t,x)|}{\eps+\eD^{-\theta |x|}}\,\leq C\,\eD^{-\eps\,\omega_\infty t}
\,\times\Big(\ &\sup_{0\leq s\leq T}\eD^{\eps\,\omega_\infty s}(|\psi'(s)|
+\|v(s,\cdot)\|_{L^\infty(\R)}
+\|\d_xv(s,\cdot)\|_{L^\infty([-x_*,x_*])})\\
&+\|(\eps+\eD^{-\theta |\cdot|})^{-1}\d_xv(0,\cdot)\|_{L^\infty(\R)}\Big)\,.
\end{align*}
\epr
\begin{proof}
We may argue separately to deal with bounds on $x\geq x_*$ on one hand and on $x\leq -x_*$ on the other hand, and provide details only for the former. From now on we focus on $x\geq x_*$.

We would like to apply Lemma~\ref{l:maximum} to both $A_\eps\,\d_x v$ and $-A_\eps\,\d_x v$ for a suitable weight $A_\eps$ equivalent to $(t,x)\mapsto\eD^{\eps\,\omega_\infty t}\,(\eps+\eD^{-\theta |x|})^{-1}$. Our choice is 
\[
A_\eps(t,x):=
\frac{\eD^{\omega_\infty\,\eps\,t}\eD^{\int_t^{+\infty}\eps\eD^{-\eps\,\omega_\infty s}\dD s}}{\eps+\eD^{-\theta\,|x|}}\,.
\]
Note that one has
\begin{align*}
\d_t(A_\eps \d_x v)
&+\left((f'(\uU_\eps+v)-\sigma_\eps+\psi')+2\frac{ \theta \eD^{-\theta |x|}}{\eps +\eD^{-\theta |x|}}\right)\d_x(A_\eps\d_xv)-\d_x^2(A_\eps\d_xv)\\
&=(f'(\uU_\eps)-f'(\uU_\eps+v)-\psi')A_\eps\uU_\eps'' +(f''(\uU_\eps)-f''(\uU_\eps+v))A_\eps\uU_\eps'^2 \\
&-A_\eps\d_xv\,\Big(\eps\eD^{-\eps\,\omega_\infty t}
-(f'(\uU_\eps+v)-\sigma_\eps+\psi')\frac{\theta \eD^{-\theta |x|}}{\eps +\eD^{-\theta |x|}}\\
&\qquad\qquad-\eps\,(\omega_\infty+g'(\uU_\eps+v))
-\frac{\theta^2\eD^{-\theta x}\,(2\,\eps\,+\eD^{-\theta x})}{(\eps+\eD^{-\theta x})^2}
+f''(\uU_\eps+v)(2\uU_\eps'+\d_xv)
\Big)\,.
\end{align*}
Fixing first $\theta>0$ sufficiently small, then $x_*$ sufficiently large and $\delta$ and $\eps$ sufficiently small, one enforce that the term in front of $\d_x(A_\eps\d_xv)$ is bounded from above away from zero and the term in front of $A_\eps\d_xv$ is bounded from below by a multiple of $\eps\eD^{-\eps\,\omega_\infty t}+\theta \eD^{-\theta |x|}$. This is sufficient to apply Lemma~\ref{l:maximum} and derive the claimed upper bound on $x\geq x_*$.
\end{proof}

\subsection{Proof of Theorem~\ref{th:main-sl}}

Our very first task when proving Theorem~\ref{th:main-sl} is to convert classical local well-posedness yielding maximal solutions $u$ to \eqref{eq:scalar-sl} into convenient local existence results for $(v,\psi')$.

This follows from the following simple observation. By design, $\Spsi(t)\equiv0$ when $0\leq t\leq 1$. Thus if $u$ solves \eqref{eq:scalar-sl} on $[0,T]\times\R$ then $(v,\psi')$ satisfying \eqref{eq:stab-ansatz}-\eqref{eq:v}-\eqref{eq:psi} may be obtained recursively through 
\begin{align*}
\psi(t)&=\psi_0\,,&
v(t,\cdot)&=u(t,\cdot-\psi(t))-\uU_\eps\,,&
\qquad\textrm{when }0\leq t\leq \min(\{1,T\})\,,&
\end{align*}
and, for any $n\in\N$,
\begin{align*}
\psi'(t)&\,=\,
\d_t\Spsi(t)\,(v_0)\,+\,\int_0^{t-1}\d_t\Spsi(t-s)\,\cN_\eps[v(s,\cdot)),\psi'(s)]\,\dD s\,,\\
\psi(t)&=\psi_0+\int_0^t\,\psi'(s)\,\dD s\,,\\
v(t,\cdot)&=u(t,\cdot-\psi(t))-\uU_\eps\,,
\qquad\textrm{when }\min(\{n,T\})\leq t\leq \min(\{n+1,T\})\,.&
\end{align*}

Now, combining together Propositions~\ref{p:Duhamel} and~\ref{p:energy}, one obtains that for any $\theta>0$ sufficiently small, there exist $\eps_0>0$, $\delta>0$ and $C\geq 1$ such that for any $0<\eps\leq \eps_0$, and $(v_0,\psi_0)$ with
\[
\|v_0\|_{L^\infty(\R)}
+\|(\eps+\eD^{-\theta |\cdot|})^{-1}\d_xv_0\|_{L^\infty(\R)}
\,\leq\,\delta\,,
\]
the corresponding solution $u$ to \eqref{eq:scalar-sl}, in the form \eqref{eq:stab-ansatz}, satisfies that if for some $T>0$ and any $0\leq t\leq T$
\begin{align*}
|\psi'(t)|
+\|v(t,\cdot)\|_{L^\infty(\R)}
&+\|(\eps+\eD^{-\theta |\cdot|})^{-1}\d_xv(t,\cdot)\|_{L^\infty(\R)}\\
&\,\leq\,2\,C\,\eD^{-\eps\,\omega_\infty t}\,
\left(\|v_0\|_{L^\infty(\R)}
+\|(\eps+\eD^{-\theta |\cdot|})^{-1}\d_xv_0\|_{L^\infty(\R)}\right)
\end{align*}
then for any $0\leq t\leq T$
\begin{align*}
|\psi'(t)|
+\|v(t,\cdot)\|_{L^\infty(\R)}
&+\|(\eps+\eD^{-\theta |\cdot|})^{-1}\d_xv(t,\cdot)\|_{L^\infty(\R)}\\
&\,\leq\,C\,\eD^{-\eps\,\omega_\infty t}\,
\left(\|v_0\|_{L^\infty(\R)}
+\|(\eps+\eD^{-\theta |\cdot|})^{-1}\d_xv_0\|_{L^\infty(\R)}\right)\,.
\end{align*}
From this and a continuity argument stem that $u$ is global and that the latter estimate holds globally in time. One achieves the proof of Theorem~\ref{th:main-sl} by deriving bounds on $\psi$ by integration of those on $\psi'$ and going back to original variables.

\appendix

\section{Wave profiles}\label{s:profiles}

In the present Appendix, we prove Proposition~\ref{p:profile}. Let us first reformulate the wave profile equation in terms of 
\begin{align*}
\tuU_\eps&:=\frac{\uU_\eps-\uU_0}{\eps}\,,&
\tsigma_\eps&:=\frac{\sigma_\eps-\sigma_0}{\eps}\,.&
\end{align*}
The equation to consider is 
\begin{align*}
\tuU_\eps''-\left((f'(\uU_0)-\sigma_0)\,\tuU_\eps\right)'
&\,=\,-\,g(\uU_0+\eps\,\tuU_\eps)
-\left(\tsigma_\eps\,(\uU_0+\eps\,\tuU_\eps)\right)'\\
&\quad+\left(\frac{f(\uU_0+\eps\,\tuU_\eps)-f(\uU_0)-\eps\,f'(\uU_0)\,\tuU_\eps}{\eps}\right)'\,,
\end{align*}
with $\tuU_\eps(0)=0$, $(\tsigma_\eps,\beD^{\theta\,|\,\cdot\,|}\tuU_\eps,\beD^{\theta\,|\,\cdot\,|}\tuU_\eps')$ uniformly bounded, for some sufficiently small $\theta>0$. As announced in the introduction the framework we first consider is suboptimal from the point of view of spatial localization but we shall refine it in a second step. To carry out the first step we introduce spaces $W^{k,\infty}_{\theta}$ and their subspaces $BUC^{k}_{\theta}$, corresponding to norms
\begin{align*}
\|v\|_{W_{\theta}^{k,\infty}(\R)}
&=\sum_{j=0}^k\,\|\eD^{\theta\,|\,\cdot\,|}\,\d_x^jv\|_{L^{\infty}(\R)}\,.
\end{align*}
In this first step, we just pick some $0<\theta<\min(\{\theta_0^\ell,\theta_0^r\})$ and let all the constants depend on this particular choice. 

We begin with two preliminary remarks. Firstly note that a simple integration yields that a necessary constraint is
\[
\tsigma_\eps\,=\,-\frac{1}{\uu_{+\infty}-\uu_{-\infty}}\,\int_{\R}g(\uU_0+\eps\,\tuU_\eps)\,=:\,\tSigma_\eps[\tuU_\eps]
\]
and that 
\begin{align*}
\tcN_\eps[\tuU_\eps]
&:=-\,g(\uU_0+\eps\,\tuU_\eps)
-\left(\tSigma_\eps[\tuU_\eps]\,(\uU_0+\eps\,\tuU_\eps)\right)'\\
&\quad+\left(\frac{f(\uU_0+\eps\,\tuU_\eps)-f(\uU_0)-\eps\,f'(\uU_0)\,\tuU_\eps}{\eps}\right)'
\end{align*}
defines a continuous map from $BUC^{1}_{\theta}$ to the closed subspace of $BUC^{0}_{\theta}$ whose range is contained in the set of functions with zero integral and that, on any ball of $BUC^{1}_{\theta}$, has an $\cO(\eps)$-Lipschitz constant.

Secondly, denoting $L_0$ the operator defined by
\begin{align*}
L_0(v)\,:=\,v''-(f'(\,\uU_0)-\sigma_0)v)'
\end{align*}
on $BUC^{0}_{\theta}$, with domain $BUC^{2}_{\theta}$, we observe that $L_0$ is Fredholm of index $0$ (as a continuous operator from $BUC^{2}_{\theta}$ to $BUC^{0}_{\theta}$), its kernel is spanned by $\uU'_\eps$ and the kernel of its adjoint is reduced to constant functions. The foregoing claims are easily proved by direct inspection but may also be obtained with the arguments of Sections~\ref{s:spectral} and~\ref{s:Green}, combining spatial dynamics point of view with a Sturm-Liouville argument. Since evaluation at $0$ acts continuously on $BUC^{2}_{\theta}$ and $\uU_0'(0)\neq0$, this implies that the restriction of $L_0$ from the closed subspace of $BUC^{2}_{\theta}$ consisting of functions with value $0$ at $0$ to the closed subspace of $BUC^{0}_{\theta}$ consisting of functions with zero integral is boundedly invertible. Indeed, the inverse of this restriction is readily seen to be given by
\begin{align*}
L_0^\dagger(h)(x):=-\int_0^x\int_{z}^{+\infty}\frac{\uU_0'(x)}{\uU_0'(z)}\ h(y)\,\dD y\,\dD z
\,=\,\int_0^x\int^{z}_{-\infty}\frac{\uU_0'(x)}{\uU_0'(z)}\ h(y)\,\dD y\,\dD z\,.
\end{align*}

Note that from the profile equation itself stems that if $\tuU_\eps$ is a $BUC^{1}_{\theta}$-solution it is also a $BUC^{2}_{\theta}$-solution so that the problem reduces to
\begin{align*}
\tsigma_\eps&=\tSigma[\tuU_\eps]\,,&
\tuU_\eps&=L_0^\dagger(\tcN_\eps[\tuU_\eps])\,.
\end{align*}
If $C_0$ is chosen such that $C_0>\|L_0^\dagger\tcN_\eps[\Zero_{\R}]\|_{W^{1,\infty}_{\theta}(\R)}$, it follows that, when $\eps$ is sufficiently small, the map $L_0^\dagger\circ\tcN_\eps$ sends the complete space
\[
\left\{\,v\in BUC^{1}_{\theta}(\R)\,;\,v(0)=0\ \textrm{and}\ \|v\|_{W^{1,\infty}_{\theta}(\R)}\leq C_0\,\right\}
\]
into itself and is strictly contracting with an $\cO(\eps)$-Lipschitz constant. Thus resorting to the Banach fixed-point theorem achieves the first step of the proof of Proposition~\ref{p:profile}.

Note that, in order to conclude the proof, it is sufficient to provide asymptotic descriptions of $\eD^{\theta_\eps^\ell|\,\cdot\,|}(\uU_\eps-\uu_{-\infty})$, $\eD^{\theta_\eps^r\,\cdot\,}(\uU_\eps-\uu_{+\infty})$, $\eD^{\theta_\eps^\ell|\,\cdot\,|}\uU_\eps'$ and $\eD^{\theta_\eps^r\,\cdot\,}\uU_\eps'$. Indeed, on one hand, the asymptotic comparisons for $\eD^{\theta_\eps^\ell|\,\cdot\,|}\uU_\eps^{(k)}$ and $\eD^{\theta_\eps^r\,\cdot\,}\uU_\eps^{(k)}$, $k\geq2$, are then deduced recursively by using the profile equation (differentiated $(k-2)$ times). On the other hand, since, for $\#\in\{\ell,r\}$, $\theta_\eps^\#=\theta_0^\#+\cO(\eps)$, the asymptotic descriptions are sufficient to upgrade the existence part of the first step arbitrarily close to optimal spatial decay rates, $\alpha^\#\to\theta_0^\#$.

As a further reduction, we observe that the asymptotics for $\eD^{\theta_\eps^\ell|\,\cdot\,|}(\uU_\eps-\uu_{-\infty})$,  $\eD^{\theta_\eps^r\,\cdot\,}(\uU_\eps-\uu_{+\infty})$, may be deduced from the ones for $\eD^{\theta_\eps^\ell|\,\cdot\,|}\uU_\eps'$ and $\eD^{\theta_\eps^r\,\cdot\,}\uU_\eps'$ by integration since
\begin{align*}
\eD^{-\theta_\eps^\ell\,x\,}
(\uU_\eps(x)-\uu_{-\infty})
&-\eD^{-\theta_0^\ell\,x\,}(\uU_0(x)-\uu_{-\infty})\\
&=\int_{-\infty}^x\eD^{\theta_\eps^\ell\,(y-x)\,}
\left(\eD^{-\theta_\eps^\ell\,y\,}\uU_\eps'(y)-\eD^{-\theta_0^\ell\,y\,}\uU_0'(y)\right)\,\dD y\\
&\quad+\int_{-\infty}^x\eD^{\theta_0^\ell\,(y-x)\,}
\left(\eD^{(\theta_\eps^\ell-\theta_0^\ell)\,(y-x)\,}-1\right)\,\eD^{-\theta_0^\ell\,y\,}\uU_0'(y)\,\dD y
\end{align*}
and likewise near $+\infty$.

To conclude, we derive the study of $\eD^{\theta_\eps^\ell|\,\cdot\,|}\uU_\eps'$ and $\eD^{\theta_\eps^r\,\cdot\,}\uU_\eps'$ from the analysis of Proposition~\ref{p:gap-lemma} (with $K=\{0\}$). Indeed, 
\begin{align*}
\theta_\eps^\ell&=\mu_+^\eps(0;\uu_{-\infty})\,,&
\theta_\eps^\ell&=-\mu_-^\eps(0;\uu_{+\infty})\,,
\end{align*}
and
\begin{align*}
\uU_\eps'(x)&=\frac{\left(\uu_{+\infty}-\uu_{-\infty}\right)}{2}\,\frac{\eD^{\theta_\eps^\ell\,x\,}\,\beD_1\cdot P^\ell_\eps(0,x)\,\bR_+^\eps(0;\uu_{-\infty})}{
\int_{-\infty}^0 \eD^{\theta_\eps^\ell\,y\,}\,\beD_1\cdot P^\ell_\eps(0,y)\,\bR_+^\eps(0;\uu_{-\infty})\,\dD y}\,,\\
&=\frac{\left(\uu_{+\infty}-\uu_{-\infty}\right)}{2}\,\frac{-\eD^{\theta_\eps^r\,x\,}\,\beD_1\cdot P^r_\eps(0,x)\,\bR_-^\eps(0;\uu_{+\infty})}{
\int_0^{+\infty} \eD^{-\theta_\eps^r\,y\,}\,\beD_1\cdot P^r_\eps(0,y)\,\bR_-^\eps(0;\uu_{+\infty})\,\dD y}\,.
\end{align*}
Thus the claimed expansion stems from the smoothness in $\eps$ afforded by Proposition~\ref{p:gap-lemma}.


\begin{thebibliography}{JNR{\etalchar{+}}19}

\bibitem[BJRZ11]{BJRZ}
B.~Barker, M.~A. Johnson, L.~M. Rodrigues, and K.~Zumbrun.
\newblock Metastability of solitary roll wave solutions of the {S}t. {V}enant
  equations with viscosity.
\newblock {\em Phys. D}, 240(16):1289--1310, 2011.

\bibitem[BGM17]{Bedrossian-Germain-Masmoudi}
J.~Bedrossian, P.~Germain, and N.~Masmoudi.
\newblock On the stability threshold for the 3{D} {C}ouette flow in {S}obolev
  regularity.
\newblock {\em Ann. of Math. (2)}, 185(2):541--608, 2017.

\bibitem[BM15]{Bedrossian-Masmoudi}
J.~Bedrossian and N.~Masmoudi.
\newblock Inviscid damping and the asymptotic stability of planar shear flows
  in the 2{D} {E}uler equations.
\newblock {\em Publ. Math. Inst. Hautes \'{E}tudes Sci.}, 122:195--300, 2015.

\bibitem[BMV16]{Bedrossian-Masmoudi-Vicol}
J.~Bedrossian, N.~Masmoudi, and V.~Vicol.
\newblock Enhanced dissipation and inviscid damping in the inviscid limit of
  the {N}avier-{S}tokes equations near the two dimensional {C}ouette flow.
\newblock {\em Arch. Ration. Mech. Anal.}, 219(3):1087--1159, 2016.

\bibitem[BGS07]{Benzoni-Serre_book}
S.~Benzoni-Gavage and D.~Serre.
\newblock {\em Multidimensional hyperbolic partial differential equations}.
\newblock Oxford Mathematical Monographs. The Clarendon Press, Oxford
  University Press, Oxford, 2007.
\newblock First-order systems and applications.

\bibitem[BB05]{Bianchini-Bressan}
S.~Bianchini and A.~Bressan.
\newblock Vanishing viscosity solutions of nonlinear hyperbolic systems.
\newblock {\em Ann. of Math. (2)}, 161(1):223--342, 2005.

\bibitem[Bre00]{Bressan}
A.~Bressan.
\newblock {\em Hyperbolic systems of conservation laws}, volume~20 of {\em
  Oxford Lecture Series in Mathematics and its Applications}.
\newblock Oxford University Press, Oxford, 2000.
\newblock The one-dimensional Cauchy problem.

\bibitem[Cro10]{Crooks}
E.~C.~M. Crooks.
\newblock Front profiles in the vanishing-diffusion limit for monostable
  reaction-diffusion-convection equations.
\newblock {\em Differential Integral Equations}, 23(5-6):495--512, 2010.

\bibitem[CM07]{Crooks-Mascia}
E.~C.~M. Crooks and C.~Mascia.
\newblock Front speeds in the vanishing diffusion limit for
  reaction-diffusion-convection equations.
\newblock {\em Differential Integral Equations}, 20(5):499--514, 2007.

\bibitem[Dav07]{Davies}
E.~B. Davies.
\newblock {\em Linear operators and their spectra}, volume 106 of {\em
  Cambridge Studies in Advanced Mathematics}.
\newblock Cambridge University Press, Cambridge, 2007.

\bibitem[DR20]{DR1}
V.~Duch{\^ e}ne and L.~M. Rodrigues.
\newblock Large-time asymptotic stability of {R}iemann shocks of scalar balance
  laws.
\newblock {\em SIAM J. Math. Anal.}, 52(1):792--820 889, 2020.

\bibitem[DRar]{DR2}
V.~Duch{\^ e}ne and L.~M. Rodrigues.
\newblock Stability and instability in scalar balance laws: fronts and periodic
  waves.
\newblock {\em Anal. PDE}, to appear.

\bibitem[Gil10]{Gilding}
B.~H. Gilding.
\newblock On front speeds in the vanishing diffusion limit for
  reaction-convection-diffusion equations.
\newblock {\em Differential Integral Equations}, 23(5-6):445--450, 2010.

\bibitem[Goo86]{Goodman}
J.~Goodman.
\newblock Nonlinear asymptotic stability of viscous shock profiles for
  conservation laws.
\newblock {\em Arch. Rational Mech. Anal.}, 95(4):325--344, 1986.

\bibitem[Goo89a]{Goodman_bis}
J.~Goodman.
\newblock Stability of viscous scalar shock fronts in several dimensions.
\newblock {\em Trans. Amer. Math. Soc.}, 311(2):683--695, 1989.

\bibitem[Goo89b]{Goodman_ter}
J.~Goodman.
\newblock Stability of viscous scalar shock fronts in several dimensions.
\newblock {\em Trans. Amer. Math. Soc.}, 311(2):683--695, 1989.

\bibitem[GX92]{Goodman-Xin}
J.~Goodman and Z.~P. Xin.
\newblock Viscous limits for piecewise smooth solutions to systems of
  conservation laws.
\newblock {\em Arch. Rational Mech. Anal.}, 121(3):235--265, 1992.

\bibitem[GR01]{Grenier-Rousset}
E.~Grenier and F.~Rousset.
\newblock Stability of one-dimensional boundary layers by using {G}reen's
  functions.
\newblock {\em Comm. Pure Appl. Math.}, 54(11):1343--1385, 2001.

\bibitem[H{\"{a}}r00]{Harterich_hyperbolic}
J.~H{\"{a}}rterich.
\newblock Viscous profiles for traveling waves of scalar balance laws: the
  uniformly hyperbolic case.
\newblock {\em Electron. J. Differential Equations}, pages No. 30, 22, 2000.

\bibitem[H{\"{a}}r03]{Harterich_canard}
J.~H{\"{a}}rterich.
\newblock Viscous profiles of traveling waves in scalar balance laws: the
  canard case.
\newblock {\em Methods Appl. Anal.}, 10(1):97--117, 2003.

\bibitem[Hen81]{Henry-geometric}
D.~Henry.
\newblock {\em Geometric theory of semilinear parabolic equations}, volume 840
  of {\em Lecture Notes in Mathematics}.
\newblock Springer-Verlag, Berlin-New York, 1981.

\bibitem[HR18]{Herda-Rodrigues}
M.~Herda and L.~M. Rodrigues.
\newblock Large-time behavior of solutions to
  {V}lasov-{P}oisson-{F}okker-{P}lanck equations: from evanescent collisions to
  diffusive limit.
\newblock {\em J. Stat. Phys.}, 170(5):895--931, 2018.

\bibitem[How99a]{Howard_lin}
P.~Howard.
\newblock Pointwise estimates on the {G}reen's function for a scalar linear
  convection-diffusion equation.
\newblock {\em J. Differential Equations}, 155(2):327--367, 1999.

\bibitem[How99b]{Howard_nonlin}
P.~Howard.
\newblock Pointwise {G}reen's function approach to stability for scalar
  conservation laws.
\newblock {\em Comm. Pure Appl. Math.}, 52(10):1295--1313, 1999.

\bibitem[HLZ09]{Humpherys-Lyng-Zumbrun}
J.~Humpherys, G.~Lyng, and K.~Zumbrun.
\newblock Spectral stability of ideal-gas shock layers.
\newblock {\em Arch. Ration. Mech. Anal.}, 194(3):1029--1079, 2009.

\bibitem[JNR{\etalchar{+}}19]{JNRYZ}
M.~A. Johnson, P.~Noble, L.~M. Rodrigues, Z.~Yang, and K.~Zumbrun.
\newblock Spectral stability of inviscid roll waves.
\newblock {\em Comm. Math. Phys.}, 367(1):265--316, 2019.

\bibitem[JNRZ14]{JNRZ-conservation}
M.~A. Johnson, P.~Noble, L.~M. Rodrigues, and K.~Zumbrun.
\newblock Behavior of periodic solutions of viscous conservation laws under
  localized and nonlocalized perturbations.
\newblock {\em Invent. Math.}, 197(1):115--213, 2014.

\bibitem[JGK93]{Jones-Gardner-Kapitula}
C.~K. R.~T. Jones, R.~Gardner, and T.~Kapitula.
\newblock Stability of travelling waves for nonconvex scalar viscous
  conservation laws.
\newblock {\em Comm. Pure Appl. Math.}, 46(4):505--526, 1993.

\bibitem[KV21a]{Kang-Vasseur_bis}
M.-J. Kang and A.~Vasseur.
\newblock Contraction property for large perturbations of shocks of the
  barotropic {N}avier-{S}tokes system.
\newblock {\em J. Eur. Math. Soc. (JEMS)}, 23(2):585--638, 2021.

\bibitem[KV21b]{Kang-Vasseur}
M.-J. Kang and A.~F. Vasseur.
\newblock Uniqueness and stability of entropy shocks to the isentropic {E}uler
  system in a class of inviscid limits from a large family of {N}avier-{S}tokes
  systems.
\newblock {\em Invent. Math.}, 224(1):55--146, 2021.

\bibitem[Kap94]{Kapitula}
T.~Kapitula.
\newblock On the stability of travelling waves in weighted {$L^\infty$} spaces.
\newblock {\em J. Differential Equations}, 112(1):179--215, 1994.

\bibitem[KP13]{KapitulaPromislow-stability}
T.~Kapitula and K.~Promislow.
\newblock {\em Spectral and dynamical stability of nonlinear waves}, volume 185
  of {\em Applied Mathematical Sciences}.
\newblock Springer, New York, 2013.
\newblock With a foreword by Christopher K. R. T. Jones.

\bibitem[Kat76]{Kato}
T.~Kato.
\newblock {\em Perturbation theory for linear operators}.
\newblock Springer-Verlag, Berlin, second edition, 1976.
\newblock Grundlehren der Mathematischen Wissenschaften, Band 132.

\bibitem[KK98]{Kreiss-Kreiss}
G.~Kreiss and H.-O. Kreiss.
\newblock Stability of systems of viscous conservation laws.
\newblock {\em Comm. Pure Appl. Math.}, 51(11-12):1397--1424, 1998.

\bibitem[Kru70]{Kruzhkov}
S.~N. Kru{\v z}kov.
\newblock First order quasilinear equations with several independent variables.
\newblock {\em Mat. Sb. (N.S.)}, 81 (123):228--255, 1970.

\bibitem[Liu85]{Liu}
T.-P. Liu.
\newblock Nonlinear stability of shock waves for viscous conservation laws.
\newblock {\em Mem. Amer. Math. Soc.}, 56(328):v+108, 1985.

\bibitem[Maj83a]{Majda1}
A.~Majda.
\newblock The existence of multidimensional shock fronts.
\newblock {\em Mem. Amer. Math. Soc.}, 43(281):v+93, 1983.

\bibitem[Maj83b]{Majda2}
A.~Majda.
\newblock The stability of multidimensional shock fronts.
\newblock {\em Mem. Amer. Math. Soc.}, 41(275):iv+95, 1983.

\bibitem[MZ03]{Mascia-Zumbrun_hyp-par}
C.~Mascia and K.~Zumbrun.
\newblock Pointwise {G}reen function bounds for shock profiles of systems with
  real viscosity.
\newblock {\em Arch. Ration. Mech. Anal.}, 169(3):177--263, 2003.

\bibitem[MZ04]{Mascia-Zumbrun_hyp-par_bis}
C.~Mascia and K.~Zumbrun.
\newblock Stability of large-amplitude viscous shock profiles of
  hyperbolic-parabolic systems.
\newblock {\em Arch. Ration. Mech. Anal.}, 172(1):93--131, 2004.

\bibitem[MN85]{Matsumura-Nishihara}
A.~Matsumura and K.~Nishihara.
\newblock On the stability of travelling wave solutions of a one-dimensional
  model system for compressible viscous gas.
\newblock {\em Japan J. Appl. Math.}, 2(1):17--25, 1985.

\bibitem[M{\'{e}}t01]{Metivier_cours-chocs}
G.~M{\'{e}}tivier.
\newblock Stability of multidimensional shocks.
\newblock In {\em Advances in the theory of shock waves}, volume~47 of {\em
  Progr. Nonlinear Differential Equations Appl.}, pages 25--103. Birkh\"{a}user
  Boston, Boston, MA, 2001.

\bibitem[MZ05]{MZ_AMS}
G.~M{\'{e}}tivier and K.~Zumbrun.
\newblock Large viscous boundary layers for noncharacteristic nonlinear
  hyperbolic problems.
\newblock {\em Mem. Amer. Math. Soc.}, 175(826):vi+107, 2005.

\bibitem[Rod13]{R}
L.~M. Rodrigues.
\newblock {\em Asymptotic stability and modulation of periodic wavetrains,
  general theory \& applications to thin film flows}.
\newblock Habilitation {\`a} diriger des recherches, Universit\'e Lyon 1, 2013.

\bibitem[Rod15]{R_Roscoff}
L.~M. Rodrigues.
\newblock Space-modulated stability and averaged dynamics.
\newblock {\em Journ\'ees \'Equations aux d\'eriv\'ees partielles},
  2015(8):1--15, 2015.

\bibitem[Rod18]{R_linKdV}
L.~M. Rodrigues.
\newblock Linear asymptotic stability and modulation behavior near periodic
  waves of the {K}orteweg--de {V}ries equation.
\newblock {\em J. Funct. Anal.}, 274(9):2553--2605, 2018.

\bibitem[RZ16]{Rodrigues-Zumbrun}
L.~M. Rodrigues and K.~Zumbrun.
\newblock Periodic-coefficient damping estimates, and stability of
  large-amplitude roll waves in inclined thin film flow.
\newblock {\em SIAM J. Math. Anal.}, 48(1):268--280, 2016.

\bibitem[Rou02]{Rousset}
F.~Rousset.
\newblock Viscous limits for strong shocks of one-dimensional systems of
  conservation laws.
\newblock In {\em Journ\'{e}es ``\'{E}quations aux {D}\'{e}riv\'{e}es
  {P}artielles'' ({F}orges-les-{E}aux, 2002)}, pages Exp. No. XVI, 12. Univ.
  Nantes, Nantes, 2002.

\bibitem[San02]{Sandstede}
B.~Sandstede.
\newblock Stability of travelling waves.
\newblock In {\em Handbook of dynamical systems, {V}ol. 2}, pages 983--1055.
  North-Holland, Amsterdam, 2002.

\bibitem[Sat73]{Sattinger-book}
D.~H. Sattinger.
\newblock {\em Topics in stability and bifurcation theory}.
\newblock Lecture Notes in Mathematics, Vol. 309. Springer-Verlag, Berlin-New
  York, 1973.

\bibitem[Sat76]{Sattinger_one}
D.~H. Sattinger.
\newblock On the stability of waves of nonlinear parabolic systems.
\newblock {\em Advances in Math.}, 22(3):312--355, 1976.

\bibitem[Sat77]{Sattinger_two}
D.~H. Sattinger.
\newblock Weighted norms for the stability of traveling waves.
\newblock {\em J. Differential Equations}, 25(1):130--144, 1977.

\bibitem[Ser21]{Serre_scalar}
D.~Serre.
\newblock Asymptotic stability of scalar multi-{D} inviscid shock waves.
\newblock {\em arXiv preprint arXiv:2103.09615}, 2021.

\bibitem[SYZ20]{SYZ}
A.~Sukhtayev, Z.~Yang, and K.~Zumbrun.
\newblock Spectral stability of hydraulic shock profiles.
\newblock {\em Phys. D}, 405:132360, 9, 2020.

\bibitem[WX05]{Wu-Xing}
Y.~Wu and X.~Xing.
\newblock The stability of travelling fronts for general scalar viscous balance
  law.
\newblock {\em J. Math. Anal. Appl.}, 305(2):698--711, 2005.

\bibitem[Xin05]{Xing}
X.-x. Xing.
\newblock Existence and stability of viscous shock waves for non-convex viscous
  balance law.
\newblock {\em Adv. Math. (China)}, 34(1):43--53, 2005.

\bibitem[YZ20]{YangZumbrun20}
Z.~Yang and K.~Zumbrun.
\newblock Stability of {H}ydraulic {S}hock {P}rofiles.
\newblock {\em Arch. Ration. Mech. Anal.}, 235(1):195--285, 2020.

\bibitem[Zum01]{Zumbrun}
K.~Zumbrun.
\newblock Multidimensional stability of planar viscous shock waves.
\newblock In {\em Advances in the theory of shock waves}, volume~47 of {\em
  Progr. Nonlinear Differential Equations Appl.}, pages 307--516.
  Birkh\"{a}user Boston, Boston, MA, 2001.

\bibitem[ZH98]{Zumbrun-Howard}
K.~Zumbrun and P.~Howard.
\newblock Pointwise semigroup methods and stability of viscous shock waves.
\newblock {\em Indiana Univ. Math. J.}, 47(3):741--871, 1998.

\bibitem[ZH02]{Zumbrun-Howard_erratum}
K.~Zumbrun and P.~Howard.
\newblock Errata to: ``{P}ointwise semigroup methods, and stability of viscous
  shock waves'' [{I}ndiana {U}niv. {M}ath. {J}. {\bf 47} (1998), no. 3,
  741--871; {MR}1665788 (99m:35157)].
\newblock {\em Indiana Univ. Math. J.}, 51(4):1017--1021, 2002.

\end{thebibliography}

\newcommand{\etalchar}[1]{$^{#1}$}
\newcommand{\SortNoop}[1]{}

\end{document}